\documentclass[a4paper,11pt]{amsart}
\usepackage{amssymb, amsmath,latexsym,amsfonts,amsbsy, amsthm,mathtools,graphicx,CJKutf8,CJKnumb,CJKulem,color,geometry}
\usepackage{verbatim}
\usepackage{mathrsfs}
\usepackage{ulem}
\usepackage{indentfirst}
\usepackage{txfonts}

\def\W{W_{n+1,\mathrm{loc}}^{2,1}}
\def\L{\mathcal{L}_\phi}
\def\I{\mathcal{I}}
\def\fz{\infty}
\def\C{C^{1+\alpha,\f{1+\alpha}{2}}_{\mathrm{loc}}}
\def\C'{C^{1+\alpha',\f{1+\alpha'}{2}}_{\mathrm{loc}}}
\def\Lfz{L^\fz}
\def\l{\left}
\def\r{\right}
\def\tz{\theta}
\def\rz{\rho}
\def\dz{\delta}
\def\bz{\beta}

\def\az{\alpha}
\def\lz{\lambda}
\def\Lz{\Lambda}
\def\az{\alpha}
\def\ez{\epsilon}
\def\gz{\gamma}

\def\f{\frac}
\def\R{\mathbb{R}^n}
\author{Lin Tang and Qian Zhang }
\title{\bf{Global $C^{1+\alpha,\f{1+\alpha}{2}}$ regularity on the linearized parabolic Monge-Amp$\grave{e}$re equation}\footnotetext{ \hspace{-0.65
cm} 2010 Mathematics Subject  Classification: 35J55.\\
The  research was supported  by the NNSF  (11271024)
and (11571289)  of China.}}
\date{}

\theoremstyle{plain}\newtheorem{thm}{Theorem}[section]
\theoremstyle{plain}\newtheorem{prop}{Proposition}[section]
\theoremstyle{plain}\newtheorem{cor}{Corollary}[section]
\theoremstyle{plain}\newtheorem{lem}{Lemma}[section]
\theoremstyle{plain}\newtheorem{rem}{Remark}[section]
\newtheorem*{thm1}{Theorem 1}
\newtheorem*{thm2}{Theorem 2}

\numberwithin{equation}{section}

\begin{document}
\maketitle
\noindent {\bf{Abstract}.}\quad
In this paper, we establish global $C^{1+\alpha,\f{1+\alpha}{2}}$ estimates for solutions of the linearized parabolic Monge-Amp$\grave{e}$re equation
$$\L u(x,t):=-u_t\,\mathrm{det}D^2\phi(x)+\mathrm{tr}[\Phi(x) D^2 u]=f(x,t)$$
under appropriate conditions on the domain, Monge-Amp$\grave{e}$re measures, boundary data and $f$, where $\Phi:=\mathrm{det}(D^2\phi)(D^2\phi)^{-1}$ is the cofactor of the Hessian of $D^2\phi$.
\bigskip

\section{Introduction}\label{s1}

This paper is concerned with global regularity of solutions of the linearized parabolic Monge-Amp$\grave{e}$re equation introduced in \cite{HQ}
\begin{equation}\label{s1:1}
\L u(x,t):=-u_t\,\mathrm{det}D^2\phi(x)+\mathrm{tr}[\Phi(x) D^2 u]=f(x,t)\quad\quad\mathrm{in}\;\Omega\times (0,T],
\end{equation}
where $\Omega$ is a convex domain in $\mathbb{R}^n$, $T>0$ and $\phi\in C(\overline{\Omega})$ is a strictly convex function satisfying
\begin{equation}\label{phi}
\mathrm{det} D^2\phi=g,\quad\quad 0<\lambda\leq g\leq\Lambda\quad\mathrm{in}\;\Omega
\end{equation}
for some constants $0<\lz\le\Lz<\infty$, $\Phi:=\mathrm{det}(D^2\phi)(D^2\phi)^{-1}$ is the cofactor of the Hessian of $D^2\phi$, $D^2 u$ denotes the Hessian of $u$ in the $x$ variable, and $\mathrm{tr}(A)$ means the trace of the matrix $A$. 

Concerning the regularity of the linearized elliptic Monge-Amp$\grave{e}$re equation
\begin{equation}\label{s1:4}
\mathrm{tr}[\Phi D^2 u]=f\quad\quad\mathrm{in}\;\Omega.
\end{equation}
a fundamental result is the Harnack inequality \cite{CG} for nonnegative solutions of 
\begin{equation*}
\mathrm{tr}[\Phi D^2 u]=0\quad\quad\mathrm{in}\;\Omega,
\end{equation*}
which yields interior H$\ddot{o}$lder continuity of solutions of \eqref{s1:4}. By using this result and perturbation arguments, Guti$\acute{e}$rrez and Nguyen \cite{GN} established interior $C^{1,\alpha}$ estimates for solutions of \eqref{s1:4}. Using the Localization Theorem \cite{S2}, global $C^{1,\az}$ regularity of solutions of the Dirichlet boundary problem
\begin{equation}\label{s1:6}
\l\{\begin{array}{rclcc}
\mathrm{tr}[\Phi D^2 u]&=&f&&\mathrm{in}\;\Omega,\\
u&=&\varphi&&\mathrm{on}\;\partial\Omega,
\end{array}
\r.
\end{equation}
was established in \cite{S1} under further assumptions on the geometry of $\Omega$ and $\phi$. 

\begin{comment}
The equation \eqref{s1:4} appears in several contexts including affine geometry \cite{TW, TW1, TW2}, complex geometry \cite{Ab, Do, Do1} and fluid mechanics \cite{B,Cu,Lo}. In particular, the authors in \cite{TW1} resolved Chern's conjecture in affine geometry concerning affine maximal hypersurfaces in $\mathbb{R}^3$.
\end{comment}

Regularity of the linearized parabolic Monge-Amp$\grave{e}$re equation \eqref{s1:1} was first studied by Huang \cite{HQ}, where a Harnack inequality for nonnegative solutions of $\L u=0$ was established on parabolic sections and thus generalizes the result in \cite{CG} from the elliptic to the parabolic case.

Hence, we are interested in developing higher regularity of solutions of the linearized parabolic Monge-Amp$\grave{e}$re equation \eqref{s1:1}. As stated in \cite{GN}, to obtain $C^{2+\alpha,1+\f{\alpha}{2}}_{\mathrm{loc}}$ estimates for the solution $u$, it is natural to assume that $\mathrm{det}\,D^2\phi$ is locally H$\ddot{o}$lder continuous. However under this condition, the equation \eqref{s1:1} becomes uniformly parabolic. Therefore, we consider the $C^{1+\alpha,\f{1+\alpha}{2}}$ regularity of solutions of \eqref{s1:1}. Precisely, we extend the results in \cite{GN} and \cite{S1} to the parabolic case. 

Our first result concerns interior $C^{1+\alpha,\f{1+\alpha}{2}}$ regularity of solutions of \eqref{s1:1}, which generalizes \cite{GN} from the elliptic to the parabolic case.

\begin{thm1}
Let $B_1\subset\Omega\subset B_n$ be a normalized convex domain, $T>0,\;Q=\Omega\times(-T,0]$. Let $\phi\in C(\overline{\Omega})$ be a convex function satisfying \eqref{phi} with $\phi=0$ on $\partial\Omega$, where $g\in C(\Omega)$. Assume that $u$ is a classical solution of $\L u=f$ in $Q$ with $[f]^{n+1}_{\alpha,Q}<\fz$ for some $0<\alpha<1$, then for any $\alpha'\in(0,\alpha)$ and any $Q'=\Omega'\times(-T',0]$ with $\Omega'\Subset \Omega$ and $0<T'<T$ we have
\begin{equation*}
\|u\|_{C^{1+\alpha',\f{1+\alpha'}{2}}(\overline{Q'})}\leq C\{\|u\|_{\Lfz(Q)}+[f]^{n+1}_{\alpha,Q}\},
\end{equation*}
where $C$ depends only on $n,\alpha,\alpha',\lambda,\Lambda,\mathrm{dist}(\Omega',\partial\Omega),T,T'$ and the modulus of continuity of $g$. 
\end{thm1}

We also establish global $C^{1+\alpha,\f{1+\alpha}{2}}$ regularity of solutions of \eqref{s1:1}, which is stated as follows.

\begin{thm2}
Assume $\Omega\subset B_{1/\rz}$ contains an interior ball of radius $\rz$ tangent to $\partial\Omega$ at each point on $\partial\Omega$. Let $\phi\in C^{0,1}(\overline{\Omega})\cap C^2(\Omega)$ be a convex function satisfying \eqref{phi}, where $g\in C(\overline{\Omega})$. Assume further that $\phi$ separates quadratically from its tangent planes on $\partial\Omega$, namely,
\begin{equation*}
\rz|x-x_0|^2\leq\phi(x)-\phi(x_0)-\nabla\phi(x_0)\cdot (x-x_0)\leq\rz^{-1}|x-x_0|^2,\quad\quad\forall x,x_0\in\partial\Omega.
\end{equation*}
Set $D:=\Omega\times (0,T]$. Let $u\in C^{2,1}(D)\cap C(\overline{D})$ be a solution of
$$\left\{
\begin{array}{rcl}
\L u=f&&{\mathrm{in}\;D},\\
u=\varphi&&{\mathrm{on}\;\partial_p D},
\end{array}\right.$$
where $\varphi\in C^{2,1}(\overline{D})$. Then $u\in C^{1+\alpha,\f{1+\alpha}{2}}(\overline{D})$ and
\begin{eqnarray*}
\|u\|_{C^{1+\alpha,\f{1+\alpha}{2}}(\overline{D})}\leq C\{\|f\|_{\Lfz(D)}+\|\varphi\|_{C^{2,1}(\overline{D})}\},
\end{eqnarray*}
where $\alpha\in(0,1)$ depends only on $n,\lambda,\Lambda,\rz$, $C>0$ depends only on $n,\lambda,\Lambda,\rz,T$ and the modulus of continuity of $g$.
\end{thm2}

In the above theorems, $[f]^{n+1}_{\alpha,Q}$ is defined in \eqref{s3:thm4.70}. For any $\az\in(0,1)$ and $D=\Omega\times(T_1,T_2]$ with $\Omega\subset\mathbb{R}^n$ and $-\infty<T_1<T_2<\infty$, the norm $\|u\|_{C^{1+\alpha,\f{1+\alpha}{2}}(\overline{D})}$ is defined in Section \ref{s2}.

The perturbation arguments for the elliptic case \cite{GN}  easily applies to the proof of interior regularity (Theorem 1), using the parabolic H$\ddot{o}$lder inequality in \cite{HQ} instead. Hence the main part of this paper is devoted to the proof of regularity near the parabolic boundary, which is made up of three parts: a) the initial surface; b) the side of $\Omega$; c) the corner of $\Omega$.

To derive an estimate for solutions $u$ of \eqref{s1:1} near a), the straightforward way is to apply again perturbation arguments as in \cite{GN}. For this, firstly, we establish pointwise H$\ddot{o}$lder continuity of $u$ at a) with the aid of auxiliary functions. This together with the interior H$\ddot{o}$lder estimates \cite{HQ} gives H$\ddot{o}$lder continuity of $u$ near a). Since we have pointwise $C^{1+\az,\f{1+\az}{2}}$ regularity estimates at a) for solutions $h$ of the good equation $\mathcal{L}_w h=0$ where $w$ is a solution of $\mathrm{det}D^2w=1$, then by comparing $u$ and $h$, and iterating the comparison process, we obtain pointwise $C^{1+\az,\f{1+\az}{2}}$ regularity of $u$ at a) for any $\az\in(0,1)$. Since the density of the Monge-Amp$\grave{e}$re measure $g=\mathrm{det}D^2\phi$ is continuous, then $\phi$ is $C^{1,\gz}_{\mathrm{loc}}$ for any $\gz\in(0,1)$. This allows us to to establish $C^{1+\az,\f{1+\az}{2}}$ regularity of $u$ near a) from the pointwise regularity above and also obtain $C^{1+\az,\f{1+\az}{2}}$ regularity of $u$ for any $\az<1$. 

To derive an estimate near b), we first need a uniform estimate of $\f{u(\cdot,t)}{d_{\partial\Omega}}$ in a neighborhood of $\partial\Omega$. This is achieved by constructing a supersolution, which is just a straightforward modification of that in \cite{S1}. Next we need the construction of subsolution to make this uniform estimate inductive. The main difficulty of is the shift of time, for this, we employ the weak Harnack inequality for nonnegative supersolutions of $\L u=0$ established in \cite{HQ}. 

To establish an estimate near the the part c), the good geometry of $\Omega$ and $\phi$ allows us to prove pointwise $C^{1+\az,\f{1+\az}{2}}$ regularity of the solutions $u$ at c) by constructing auxiliary functions. This together with the estimates of $u$ at b) implies regularity estimates near c).

The organization of the paper is as follows. In Section \ref{s2}, some properties of sections are collected and interior H$\ddot{o}$lder estimates of the parabolic version are given. In Section \ref{s3}, interior $C^{1+\alpha,\f{1+\alpha}{2}}$ regularity for solutions of \eqref{s1:1} is established and Theorem $1$ is proved completely. In Sections \ref{s4}-\ref{s6}, we prove regularity for solutions of the Dirichlet problem in Theorem $2$ near the initial surface a), the side of $\Omega$ b), and the corner of $\Omega$ c) respectively. In Section \ref{s7}, we give the complete proof of Theorem $2$.

\section{Preliminaries}\label{s2}
For a domain $D=\Omega\times (T_1,T_2]$ with $\Omega\subset\mathbb{R}^{n}$ and $-\infty<T_1<T_2<\infty$, the parabolic boundary of $D$ is defined by
$$\partial_p D:=\l(\overline{\Omega}\times \{T_1\}\r)\bigcup\l(\partial\Omega\times (T_1,T_2)\r).$$
For a function $u\in C^{1+\az,\f{1+\az}{2}}(\overline{D})$ with $\az\in(0,1)$, we define
\begin{eqnarray*}
\|u\|_{C^{1+\alpha,\f{1+\alpha}{2}}(\overline{D})}:&=&\|u\|_{L^\infty(D)}+\sup_{\substack{(x,t),(x,s)\in D\\ t\neq s}}\f{|u(x,t)-u(x,s)|}{|t-s|^{\f{1+\alpha}{2}}}\\
&+&\|D u\|_{L^\infty(D)}+\sup_{\substack{(x,t),(y,s)\in D\\ (x,t)\neq (y,s)}}\f{|D u(x,t)-D u(y,s)|}{(|x-y|+|t-s|^{\f{1}{2}})^{\alpha}}.
\end{eqnarray*}

Let $\Omega$ be a convex domain in $\mathbb{R}^n$ and $\phi\in C(\overline{\Omega})$ be a strictly convex function. A section of $\phi$ at $x_0\in\overline{\Omega}$ with height $h$ is defined by
$$S_{\!\phi}(x_0,h):=\{x\in\overline{\Omega}:\phi(x)<\phi(x_0)+\nabla\phi(x_0)\cdot(x-x_0)+h\},$$
when $x_0\in\partial\Omega$, the term $\nabla\phi(x_0)$ is understood in the sense that
$$x_{n+1}=\phi(x_0)+\nabla\phi(x_0)\cdot(x-x_0)$$
is a supporting hyperplane for the graph of $\phi$ but for $\epsilon>0$,
$$x_{n+1}=\phi(x_0)+(\nabla\phi(x_0)+\epsilon\nu_{x_0})\cdot(x-x_0)$$
is not a supporting hyperplane, where $\nu_{x_0}$ denotes the interior unit normal to $\partial\Omega$ at $x_0$. Denote for simplicity
$$l_{\phi,x_0}:=\phi(x_0)+\nabla\phi(x_0)\cdot(x-x_0).$$

For $z_0=(x_0,t_0)$ with $x_0\in\overline{\Omega} ,\,t_0\in\mathbb{R}$, a parabolic section $Q_{\phi}(z_0,h)$ is defined by
$$Q_{\phi}(z_0,h):=S_{\!\phi}(x_0,h)\times (t_0-h,t_0].$$

The linearized parabolic Monge-Amp$\grave{e}$re operator $\mathcal{L}_\phi$ related to $\phi$ is defined in \eqref{s1:1}. 

Let $\Omega'\Subset\Omega$. Throughout this paper when we say a constant depends on $\Omega'$ we always mean that the dependence is only on $\mathrm{dist}(\Omega',\partial\Omega)$. For a finite number of constants $C,C_1,\dots,C_k$ ($k\in\mathbb{N}_+$), when we write $C=C(C_1,\dots,C_k)$ we always mean that $C$ depends only on $C_1,\dots,C_k$.

For any function $u$ defined in $B_{\epsilon}(x_0)\times\{t_0\}$ for some $\epsilon>0$ such that the gradient of $u$ in the $x$ variable $Du(x_0,t_0)$ exists, we always denote
$$l_{u,(x_0,t_0)}:=u(x_0,t_0)+Du(x_0,t_0)\cdot(x-x_0).$$

For any $r>0$, let $B_r(x)$ be the ball in $\mathbb{R}^n$ centered at $x$ with radius $r$ and $B_r:=B_r(0)$, we always denote $Q_r:=B_r\times(-r^2,0]$.

In establishing interior regularity and regularity near the initial surface (Sections \ref{s3}-\ref{s4}), we use one of the two hypotheses below frequently: 

\noindent$\mathbf{(H)}$ $B_1\subset\Omega\subset B_n$ is a normalized convex domain and $\phi\in C(\overline{\Omega})$ is a solution of \eqref{phi} with $\phi=0$ on $\partial\Omega$. 

\noindent$\mathbf{(H')}$ $B_{\f{6}{5}}\subset\Omega\subset B_n$ is a normalized convex domain and $\phi\in C(\overline{\Omega})$ is a solution of \eqref{phi} with $\phi=0$ on $\partial\Omega$. 

The following Lemmas \ref{s2:GN2}-\ref{s2:GN1} hold under the assumption $\mathbf{(H)}$ or $\mathbf{(H')}$. See \cite[Theorem 3.3.7, Corollary 3.2.4]{G}, and \cite[Lemma 2.1]{GN} respectively.
\begin{lem}\label{s2:GN2}
There exists $\tz_0=\tz_0(n,\lz,\Lz)>1$ such that if $x\in S_{\!\phi}(y,h)$ and $S_{\!\phi}(y,2 h)\Subset\Omega$, then $S_{\!\phi}(y,h)\subset S_{\!\phi}(x,\tz_0 h)$.
\end{lem}

\begin{lem}\label{s2:GN3}
There exist a constant $C=C(n,\lz,\Lz)$ such that $C^{-1}h^{n/2}\leq |S_{\!\phi}x,h)|\leq Ch^{n/2}$ whenever $S_{\!\phi}(x,h)\Subset\Omega$.
\end{lem}

\begin{lem}\label{s2:GN1}
For any $\Omega'\Subset\Omega$, there exist positive constants $h_0, C$ and $b$ such that for $x_0\in\Omega'$, and $0<h\leq h_0$,
\begin{equation*}
B_{C^{-1}h}(x_0)\subset S_{\!\phi}(x_0,h)\subset B_{Ch^b}(x_0),
\end{equation*}
where $b=b(n,\lambda,\Lambda)$, $h_0, C$ depend only on $n,\lambda,\Lambda$ and $\Omega'$.
\end{lem}

The following Aleksandrov-Bakelman-Pucci type maximum principle was established, for example, in \cite{Tso} or \cite[Theorem 7.1]{GM}.

\begin{lem}\label{s2:GM}
Let $D=\Omega\times (0,T]$ for some domain $\Omega\subset\R$. Assume $u\in\W(D)\bigcap C(\overline{D})$ satisfies $Lu:=-u_t+a^{ij}D_{ij}u\geq f$ in $D$. Then 
\begin{equation*}
\sup_D u\leq \sup_{\partial_p D}u^+ + C(n)(\mathrm{diam}(\Omega))^{\f{n}{n+1}}\|f/(\mathrm{det}\,a^{ij})^{\f{1}{n+1}}\|_{L^{n+1}(D)}.
\end{equation*}
\end{lem}

The Harnack inequality established in \cite{HQ} (see \cite[Remark 4.2]{HQ}) implies the oscillation estimate below which we formulate for the inhomogeneous equation by Lemma \ref{s2:GM}.

\begin{thm}\label{s2:H3}
Assume $\mathbf{(H)}$ holds. Let $u$ be a classical solution of $\L u=f$ in $Q_1$. Let $z_0\in Q_{\f{3}{4}}$ and $R\le h_1$. Then
\begin{equation}
osc_{Q_\phi(z_0,\rz)}u\leq C\l(\f{\rz}{R}\r)^{\alpha}\l(\sup_{Q_\phi(z_0,R)}|u|+R^{\f{n}{2(n+1)}}\|f\|_{L^{n+1}(Q_\phi(z_0,R))}\r),
\end{equation}
for all $\rz<R$, where $\alpha\in(0,1), h_1,C>0$ depend only on $n,\lambda,\Lambda$.
\end{thm}
 
Theorem \ref{s2:H3} implies the following interior H$\ddot{o}$lder estimates.
\begin{cor}\label{s2:H5}
Assume $\mathbf{(H)}$ holds. Let $u$ be a classical solution of $\L u=f$ in $Q_1$. Then for all $(x,t), (y,s)\in Q_{\f{3}{4}}$ we have
\begin{eqnarray*}
|u(x,t)-u(y,s)|&\leq&C (|x-y|+|t-s|^{1/2})^{\beta}\l[\|u\|_{\Lfz(Q_1)}+\|f\|_{L^{n+1}(Q_1)}\r],
\end{eqnarray*}
where $\beta\in(0,1),C>0$ depend only on $n,\lambda,\Lambda$.
\end{cor}

Let $\Omega$ be a normalized convex domain. Assume that $w$ is the convex solution of the equation
\begin{equation}\label{w}
\left\{
\begin{array}{rcl}
\mathrm{det} D^2 w=1&&{\mathrm{in}\;\Omega},\\
w=0&&{\mathrm{on}\;\partial\Omega},
\end{array}\right.
\end{equation}
then the operator $\mathcal{L}_w h=-h_t+\mathrm{tr}[(D^2 w)^{-1}D^2 h]$ is uniformly parabolic in the interior of $\Omega$ from Pogorelov's estimate. Hence we have the classical $C^{2,1}$ and H$\ddot{o}$lder estimates for uniformly parabolic equations below, see for example \cite[Chapter 3]{F} and \cite[Theorem 10.1]{LS}.

\begin{lem}\label{s2:S}
Let $B_{\f{6}{5}}\subset\Omega\subset B_n$ be a normalized convex domain, and $w$ be a solution of \eqref{w}. Then for any $\varphi\in C(\partial_p(Q_1))$, there exists a unique solution $h\in C^{2,1}(Q_1)\bigcap C(\overline{Q_1})$ of $\mathcal{L}_w h=0$ in $Q_1$ and $h=\varphi$ on $\partial_p(Q_1)$ such that
\begin{equation}\label{s2:S1}
\|h\|_{C^{2,1}(\overline{Q_{\f{3}{4}}})}:=\sum_{k=0}^2\|D_x^k h\|_{C(\overline{Q_{\f{3}{4}}})}+\|D_t h\|_{C(\overline{Q_{\f{3}{4}}})}\leq K\|\varphi\|_{\Lfz(\partial_p(Q_1))}.
\end{equation}
If in addition, $\varphi\in C^{\bz,\f{\bz}{2}}(\partial_p(Q_1))$ for some $0<\bz<1$, then $h\in C^{\alpha,\f{\alpha}{2}}(\overline{Q_1})$, and
\begin{eqnarray}\label{s2:S2}
\|h\|_{C^{\alpha,\f{\alpha}{2}}(\overline{Q_1})}\leq\bar{K}\|\varphi\|_{C^{\bz,\f{\bz}{2}}\l(\partial_p(Q_1)\r)},
\end{eqnarray}
where $\az\in(0,\bz),\bar{K}=\bar{K}(n,\bz), K=K(n)>0$, and
$$\|h\|_{C^{\alpha,\f{\alpha}{2}}(\overline{Q_1})}:=\|h\|_{C(\overline{Q_{1}})}+\sup_{\substack{(x,t),(y,s)\in Q_1,\\ (x,t)\neq (y,s)}}\f{|h(x,t)-h(y,s)|}{(|x-y|+|t-s|^{\f{1}{2}})^{\alpha}}.$$
\end{lem}

The following lemma concerns the eccentricity of sections of solutions $\phi$ of \eqref{phi}, and is used in establishing interior estimates and estimates near the initial surface (Sections \ref{s3}-\ref{s4}). See \cite[Lemmas 3.2,3.3]{GN}.

\begin{lem}\label{s2:GN5}
There exist constants $c_0=c_0(n),C_0=C_0(n)$ such that the following hold:
\begin{enumerate}
\item[\rm(i)] If the hypothesis $\mathbf{(H)}$ holds with $\lambda=1-\tz,\Lambda=1+\tz$, then there exists a positive definite matrix $M=A^t A$ satisfying
\begin{equation*}
\mathrm{det}\;M=1,\quad 0<c_0 I\leq M\leq C_0 I,
\end{equation*}
such that for $0<\mu\leq c_0$ and $\tz\leq c_0\mu^2$,  we have
\begin{equation*}\label{s2:GN51}
B_{[1-C_0(\mu^{\f{1}{2}}+\mu^{-1}\tz^{\f{1}{2}})]\sqrt{2}}\subset\mu^{-\f{1}{2}}T S_{\!\phi}(x_0,\mu)\subset B_{[1+C_0(\mu^{\f{1}{2}}+\mu^{-1}\tz^{\f{1}{2}})]\sqrt{2}}.
\end{equation*}
Moreover, let $w$ be a convex solution of \eqref{w}, then
\begin{equation*}\label{s2:GN53}
S_{\!\phi}(x_0,\mu)\subset S_{\!w}(x_0,\mu+C_0\theta^{\f{1}{2}})\subset B_{C_0\sqrt{\mu+C_0\theta^{1/2}}}(x_0),
\end{equation*}
where $x_0\in\Omega$ is the minimum point of $\phi,\;Tx:=A(x-x_0)$.
\item[\rm(ii)]Let $B_{(1-\sigma)\sqrt{2}}\subset\Omega\subset B_{(1+\sigma)\sqrt{2}}$ be a convex domain with $0<\sigma\leq 1/4$ and $\phi\in C(\overline{\Omega})$ be a convex solution of \eqref{phi} with $\lambda=1-\tz,\Lambda=1+\tz$. Then there exists a positive definite matrix $M=A^t A$ satisfying
\begin{equation*}
\mathrm{det}\;M=1,\quad 0<(1-C_0\sigma)I\leq M\leq (1+C_0\sigma)I,
\end{equation*}
such that for $0<\mu\leq c_0$ and $\tz\leq c_0\mu^2$,  we have
\begin{equation*}\label{s2:GN52}
B_{[1-C_0(\sigma\mu^{\f{1}{2}}+\mu^{-1}\tz^{\f{1}{2}})]\sqrt{2}}\subset\mu^{-\f{1}{2}}T S_{\!\phi}(x_0,\mu)\subset B_{[1+C_0(\sigma\mu^{\f{1}{2}}+\mu^{-1}\tz^{\f{1}{2}})]\sqrt{2}},
\end{equation*}
where $x_0\in\Omega$ is the minimum point of $\phi,\;Tx:=A(x-x_0)$.
\end{enumerate}
\end{lem}

\section{Interior regularity}\label{s3}

To prove Theorem $1$ we follow the perturbation arguments developed in \cite{CC,GN}. Firstly, by applying the results concerning the convergence of cofactor matrices in $L^p$ proved in \cite{GN}, we compare solutions $u$ of \eqref{s1:1} and $h$ of $\mathcal{L}_w h=0$ with the same boundary data on the parabolic boundary, where $w$ is the convex solution of \eqref{w}. Next, we iterate the comparison process to establish the regularity of the solution $u$ at the minimum point $z_0=(x_0,0)$ of the parabolically convex function $\phi(x)-t$. The group of affine transformations $A T(n)\times A T(1)$ is applied to parabolic sections $Q_{\phi}(z_0,h)$ centered at $z_0$ when rescaling the solutions. We apply the results concerning the eccentricity of sections of $\phi$ under affine transformations in \cite{GN} to $S_{\!\phi}(x_0,h)$ and in the $t$ direction we only perform a corresponding parabolic dilation, consequently the parabolic sections are more like the usual parabolic cylinder $Q_1=B_1\times (-1,0]$. Since the proof is very similar to the elliptic case \cite{GN}, we just list the parabolic versions of the main lemmas in \cite{GN} and we use them to sketch the proof of Theorem $1$.

In this section we always assume $\mathbf{(H)}$ or $\mathbf{(H')}$ holds and take $v,h$ to be solutions of
\begin{equation}\label{v}
\L v=f\quad\mathrm{in}\;Q_1,
\end{equation}
and 
\begin{equation}\label{h}
\left\{
\begin{array}{rcl}
\mathcal{L}_w h=0&&{\mathrm{in}\;Q_{1}}\\
h=v&&{\mathrm{on}\;\partial_p(Q_{1})},
\end{array}\right.
\end{equation}
respectively, where $w$ is a solution of \eqref{w}.

Comparing $v$ and $h$ we can obtain the following variant of \cite[Lemma 4.1]{GN}, using the interior H$\ddot{o}$lder estimate for \eqref{s1:1} (Corollary \ref{s2:H5}), the parabolic ABP maximum principle (Lemma \ref{s2:S}), the $C^{2,1}$ estimate for \eqref{h} (Lemma \ref{s2:GM}) instead.

\begin{lem}\label{s3:lem4.1}
Assume $\mathbf{(H')}$ holds with $\lambda=\f{1}{2},\Lambda=\f{3}{2}$. Let $v$ be a solution of \eqref{v} satisfying 
\begin{equation}\label{s3:rem4.11}
|v(x,t)-v(y,s)|\leq H(|x-y|+|t-s|^{\f{1}{2}})^{\alpha_0}\quad\quad\forall (x,t),(y,s)\in\overline{Q_1}
\end{equation}
for some $\alpha_0\in(0,1]$ and $H>0$. Assume $h\in C^{2,1}(Q_{1})\bigcap C(\overline{Q_{1}})$ is a solution of \eqref{h}. Then for any $\tau\in(0,1)$, we have
\begin{eqnarray*}
\|v-h\|_{\Lfz\l(B_{1-\tau}\times(-1+\tau^2,0]\r)}
&\leq&C\l(\|\Phi-W\|_{L^{n+1}(B_1)}+\|\mathrm{det}D^2\phi-\mathrm{det}D^2 w\|_{L^{n+1}(B_1)}\r)^{\alpha}\\
&\cdot&\l[\|v\|_{\Lfz(Q_1)}+H\r]+C\|f\|_{L^{n+1}(Q_1)},
\end{eqnarray*}
whenever $\|\Phi-W\|_{L^{n+1}(B_1)}+\|\mathrm{det}D^2\phi-\mathrm{det}D^2 w\|_{L^{n+1}(B_1)}\leq \tau^2$. Here $\alpha\in(0,\f{\alpha_0}{2}]$ depends only on $n,\alpha_0$; $C=C(n,H,\az_0)>0$.
\end{lem}

The lemma below follows from Lemma \ref{s3:lem4.1} and the results concerning the convergence of cofactor matrices of $D^2\phi$ in $L^p$ in \cite{GN}, using similar arguments as in \cite[Lemma 4.2]{GN}.

\begin{lem}\label{s3:lem4.2}
For any $M,H,\epsilon>0,\alpha_0\in(0,1]$, there exists $\dz=\dz(n,\epsilon,M,H,\alpha_0)>0$ such that if the hypothesis $\mathbf{(H')}$ holds with $\lambda=1-\dz,\Lambda=1+\dz$, $v$ is a solution of \eqref{v} with $\|v\|_{\Lfz(Q_1)}\leq M,\;\|f\|_{L^{n+1}(Q_1)}\leq\dz$ and
\begin{equation}\label{s3:rem4.21}
|v(x,t)-v(y,s)|\leq H(|x-y|+|t-s|^{\f{1}{2}})^{\alpha_0}\quad\quad\forall (x,t),(y,s)\in\overline{Q_1},
\end{equation}
and $h$ is a solution of \eqref{h}, then
\begin{eqnarray*}
\|v-h\|_{\Lfz\l(Q_{1}\r)}\leq\epsilon.
\end{eqnarray*}
\end{lem}

For a strictly convex function $\phi\in C(\overline{\Omega})$, assume for simplicity $0$ is the minimum point of $\phi$. For any $r>0$, we write $S_{\!r}(\phi):=S_{\!\phi}(0,r)$ and $Q_{r}(\phi):=S_{\!r}(\phi)\times (-r,0]$. Assume $\mathbf{(H)}$ holds with $\lz=1-\tz,\Lz=1+\tz$. Let $A$ be the affine transformation given by Lemma \ref{s2:GN5} (i). We define $\mathcal{I}(x,t):=(Ax,t)$, then by Lemma \ref{s2:GN5} (i)
$$B_{[1-\dz_1]\sqrt{2}}\times(-1,0]\subset\mu^{-\f{1}{2}}\mathcal{I}Q_{\mu}(\phi)\subset B_{[1+\dz_1]\sqrt{2}}\times(-1,0]$$
if $0<\mu\le c_0$ and $\tz\le c_0\mu^2$, where $\dz_1:=C_0(\mu^{\f{1}{2}}+\mu^{-1}\tz^{\f{1}{2}})$, the dilation is with respect to $(0,0)\in\mathbb{R}^{n}\times\mathbb{R}$, and for any $(x^*,t^*)\in\mathbb{R}^{n}\times\mathbb{R},\eta>0, D\subset\mathbb{R}^{n+1}$, the parabolic dilation with respect to $(x^*,t^*)$ is defined as 
\begin{equation}\label{d}
\eta D:=\{(x^*+\eta(x-x^*),t^*+\eta^2(t-t^*)): (x,t)\in D\} 
\end{equation}
Define
$$\phi^*(x):=\f{1}{\mu}[\phi(\mu^{\f{1}{2}}A^{-1}x)-\phi(0)-\mu],\quad x\in\Omega^*:=\mu^{-\f{1}{2}}AS_{\!\mu}(\phi).$$
For a solution $u$ of \eqref{s1:1}, define
$$v(x,t):=u(\mu^{\f{1}{2}}A^{-1}x,\mu t),\quad (x,t)\in\Omega^*\times(-1,0].$$
Comparing $v$ and the solution $h$ of \eqref{h} by Lemma \ref{s3:lem4.2} and iterate the comparison process we cant obtain the regularity of the solution $u$ at the minimum point $(0,0)$.

\begin{thm}\label{s3:thm4.5}
Given $0<\alpha'<\alpha<1,\; r_0,C_1^*>0,\; 0<T_1<T<T_2$. Assume the hypothesis $\mathbf{(H)}$ holds with $\lambda=1-\tz,\Lambda=1+\tz$, and $u$ is a classical solution of $\L u=f$ in $Q:=\Omega\times(-T,0]$ with
\begin{equation*}
\l(\f{1}{|Q_r(\phi)|}\int\!\!\!\int_{Q_r(\phi)}|f|^{n+1}dxdt\r)^{\f{1}{n+1}}\leq C_1^* r^{\f{\alpha-1}{2}}\quad\quad \mathrm{for\;all}\;Q_r(\phi)=S_{\!\phi}(x_0,r)\times (-r,0]\Subset Q\;\mathrm{with}\;r\leq r_0,
\end{equation*}
where $x_0$ is the minimum point of $\phi$. Then $u$ is $C^{1+\alpha',\f{1+\alpha'}{2}}$ at $(x_0,0)$, more precisely, there is an affine function $l(x)$ such that
\begin{equation*}
r^{-(1+\alpha')}\|u-l\|_{\Lfz(B_r(x_0)\times (-r^2,0])}+|Dl|\leq C[\|u\|_{\Lfz(Q)}+C_1^*]\quad\quad\forall r\leq \mu^*,
\end{equation*}
where $\tz\in(0,1), C,\mu^*>0$ depend only on $n,\alpha,\alpha',r_0,T_1,T_2$.
\end{thm}

Let $T>0,\;Q=\Omega\times(-T,0]$. For $f\in L^{n+1}_{\mathrm{loc}}(Q)$, define
\begin{equation}\label{s3:thm4.70}
[f]^{n+1}_{\alpha,Q}:=\sup_{Q_{\phi}(z,r)\Subset Q}r^{\f{1-\alpha}{2}}\l(\f{1}{|Q_{\phi}(z,r)|}\int\!\!\!\int_{Q_{\phi}(z,r)}|f|^{n+1}dxdt\r)^{\f{1}{n+1}}.
\end{equation}

\noindent$\mathbf{Proof\;of\;Theorem\;1.}$ Since $\Omega'\Subset \Omega$, by Lemma \ref{s2:GN1}, for any $\epsilon_0>0$, there exist constants $C,b,h_0>0$ depending only on $n,\lambda,\Lambda$, $\Omega',T,T'$ and the modulus of continuity of $g$ such that for each $z_0=(x_0,t_0)\in\Omega'\times (-T',0]$, \begin{equation*}
B_{C^{-1}h_0}(x_0)\subset S_{\!\phi}(x_0,h_0)\subset B_{Ch_0^b}(x_0),\quad\quad |g(y)-g(x_0)|\leq \epsilon_0\quad\forall y\in S_{\!\phi}(x_0,h_0).
\end{equation*}
Let $Tx=A(x-x_0)+y_0$ be an affine transformation such that $B_1\subset T S_{\!\phi}(x_0,h_0)\subset B_n$. Denote $\kappa_0:=\f{|\mathrm{det}A|^{\f{2}{n}}}{g(x_0)^{\f{1}{n}}}$ and $T_{\!p}(x,t):=(Tx,\kappa_0(t-t_0))$, it follows that
\begin{equation}\label{s3:thm4.72}
B_1\times (-\kappa_0\,h_0,0]\subset T_{\!p}(Q_{\phi}(z_0,h_0))
=T S_{\!\phi}(x_0,h_0)\times (-\kappa_0\,h_0,0]\subset B_n\times (-\kappa_0\,h_0,0],
\end{equation}
where $c(n,\lz,\Lz)\le\kappa_0h_0\le C(n,\lz,\Lz)$. Define
\begin{equation*}
\phi^*(y)=\kappa_0[(\phi-l_{\phi,x_0})(T^{-1}y)-h_0],
\end{equation*}
and
\begin{equation*}
v(y,s)=g(x_0)\kappa_0^{\f{1+\alpha}{2}}u(T_{\!p}^{-1}(y,s)).
\end{equation*}
Applying Theorem \ref{s3:thm4.5} to $v$ and by similar computation as in the proof of \cite[Theorem 4.7]{GN}, we can obtain the conclusion of Theorem 1.
\qed

\section{Regularity near the initial surface}\label{s4}

In this section we establish regularity near the initial surface for solutions of \eqref{s1:1}. Here we always assume $\mathbf{(H)}$ or $\mathbf{(H')}$ holds and we take $u$ to be a solution of
\begin{equation}\label{u}
\left\{
\begin{array}{rcl}
\mathcal{L}_\phi u=f&&\mathrm{in}\;D:=\Omega\times(0,T],\\
u=\varphi&&\mathrm{on}\;\Omega\times\{0\}.
\end{array}\right.
\end{equation}

Let $x_0\in\Omega$ and $h>0$. We now work on parabolic sections centered at $(x_0,h)$ and denote for simplicity $Q'_{\phi}(x_0,h):=Q_{\phi}((x_0,h),h)=S_{\!\phi}(x_0,h)\times(0,h]$.

The lemma below is a pointwise H$\ddot{o}$lder estimate for $u$ on $\Omega\times\{0\}$.

\begin{lem}\label{s4:lem1}
Assume $\mathbf{(H)}$ holds with $\lambda=\f{1}{2},\Lambda=\f{3}{2}$, $x_0\in\Omega'\Subset\Omega, h\le h_0$. Let $u$ be a continuous solution of \eqref{u} with
$$\|u\|_{\Lfz(Q'_{\phi}(x_0,h))}+\|f\|_{\Lfz(Q'_{\phi}(x_0,h))}+\|\varphi(\cdot,0)\|_{C^{0,1}(\overline{Q'_{\phi}(x_0,h)})}\leq 1.$$
Then for any $(x,t)\in Q'_{\phi}(x_0,h)$, we have
\begin{eqnarray*}
|u(x,t)-\varphi(x_0,0)|\leq Ch^{-\f{b}{1-b}}[\phi_{x_0}(x)+t]^{b},
\end{eqnarray*}
where $\phi_{x_0}:=\phi-l_{\phi,x_0}$, $C=C(n,\Omega')$, $h_0=h_0(n,\Omega',T)$ and $b=b(n)$. 

\begin{proof}
Let $v(x,t):=\phi_{x_0}(x)+\mu t$ with $\mu\geq n+2$. We have
\begin{eqnarray}\label{s4:lem12}
\L v=-\mu\mathrm{det} D^2\phi+n\mathrm{det} D^2\phi\leq -(\mu-n)/2\leq -1.
\end{eqnarray}
By Lemma \ref{s2:GN1}, there exist $h_0, C>0$ depending only on $n,\Omega'$ and $b=b(n)$ such that for any $0<h\le h_0$, we have
\begin{eqnarray}\label{s4:lem13}
v(x,0)=\phi_{x_0}(x)\geq (C^{-1}|x-x_0|)^{\f{1}{b}}\quad\quad x\in S_{\!\phi}(x_0,h).
\end{eqnarray}
For any $\epsilon>0$, by straightforward computation we find that 
$$K_1\epsilon^{1-\f{1}{b}}\geq r^{1-\f{1}{b}}-\epsilon\,r^{-\f{1}{b}},\quad\quad\forall r>0,$$
where $K_1:=(1-b)^{\f{1}{b}}\f{b}{1-b}$. 

Fix $\ez>0$ such that $K_1\ez^{1-\f{1}{b}}\ge 2h^{-1}\ge 1$. Then for any $x\in S_{\!\phi}(x_0,h)$, by \eqref{s4:lem13},
\begin{eqnarray}\label{s4:lem14}
|x-x_0|\leq\epsilon+K_1\ez^{1-\f{1}{b}}|x-x_0|^{\f{1}{b}}\leq\epsilon+K_1\ez^{1-\f{1}{b}}\,C^{\f{1}{b}}v(x,0).
\end{eqnarray}
Define
$$\omega^{\pm}(x,t):=\epsilon+K_1\ez^{1-\f{1}{b}}C^{\f{1}{b}}v(x,t)\pm[u(x,t)-\varphi(x_0,0)].$$
By \eqref{s4:lem12}, noting that $|f|\leq 1$, we have
$$\L\omega^{\pm}=K_1\ez^{1-\f{1}{b}}\,C^{\f{1}{b}}\L v\pm f\leq -C^{\f{1}{b}}+1\leq 0\quad\quad\mathrm{in}\;Q'_{\phi}(x_0,h).$$
On $S_{\!\phi}(x_0,h)\times\{0\}$, recall that $\|\varphi(\cdot,0)\|_{C^{0,1}(\overline{Q'_{\phi}(x_0,h)})}\le 1$, then the estimate \eqref{s4:lem14} implies that
$$\omega^{\pm}(x,0)=\epsilon+K_1\ez^{1-\f{1}{b}}\,C^{\f{1}{b}}v(x,0)\pm[\varphi(x,0)-\varphi(x_0,0)]\geq 0.$$
On $\partial S_{\!\phi}(x_0,h)\times (0,h]$, we use the fact that $\|u\|_{\Lfz(Q'_{\phi}(x_0,h))}\leq 1$ and find
$$\omega^{\pm}(x,t)\geq K_1\ez^{1-\f{1}{b}}\,C^{\f{1}{b}}\phi_{x_0}(x)-2\|u\|_{\Lfz(Q'_{\phi}(x_0,h))}\geq 0,$$
where we recall $K_1\ez^{1-\f{1}{b}}\ge 2h^{-1}$. Hence the maximum principle gives
$$\omega^{\pm}\geq 0\quad\quad\mathrm{in}\;Q'_{\phi}(x_0,h)$$
or
\begin{eqnarray}\label{s4:lem15}
|u(x,t)-\varphi(x_0,0)|\leq\epsilon+K_1\ez^{1-\f{1}{b}}\,C^{\f{1}{b}}v(x,t),\quad\quad\forall (x,t)\in Q'_{\phi}(x_0,h).
\end{eqnarray}

For $(x,t)$ such that
$$0<v(x,t)^b\le (2h^{-1}K_1^{-1})^{\f{b}{b-1}},$$
choose $\epsilon(x,t)=v(x,t)^{b}$, then we have by \eqref{s4:lem15} 
\begin{eqnarray*}
|u(x,t)-\varphi(x_0,0)|\leq K_2v(x,t)^{b}\leq K_2\mu[\phi_{x_0}(x)+t]^{b},
\end{eqnarray*}
where $K_2:=1+C^{\f{1}{b}}K_1$. 

On the other hand, for $(x,t)$ such that
$$v(x,t)^b\geq (2h^{-1}K_1^{-1})^{\f{b}{b-1}},$$
we have
\begin{eqnarray*}
|u(x,t)-\varphi(x_0,0)|\leq 2(2h^{-1}K_1^{-1})^{\f{b}{1-b}}v(x,t)^{b}=K_3 h^{-\f{b}{1-b}}[v(x,t)]^{b},
\end{eqnarray*}
where the constant $K_3$ depends only on $b$. Combining the two cases we get the desired result.
\end{proof}
\end{lem}

Lemma \ref{s4:lem1} together with interior H$\ddot{o}$lder estimates Corollary \ref{s2:H3} gives H$\ddot{o}$lder estimates for solutions of \eqref{u} near the initial surface.

\begin{lem}\label{s4:lem2}
Assume $\mathbf{(H)}$ holds with $\lambda=\f{1}{2},\Lambda=\f{3}{2}$, $x_0\in\Omega'\Subset\Omega$ and $0<h\leq h_0$. Let $u$ be a continuous solution of \eqref{u} with
$$\|u\|_{\Lfz(Q'_{\phi}(x_0,\tz_0^2 h))}+\|f\|_{\Lfz(Q'_{\phi}(x_0,\tz_0^2 h))}+\|\varphi(\cdot,0)\|_{C^{0,1}(\overline{Q'_{\phi}(x_0,\tz_0^2 h)})}\le 1,$$
where $\tz_0$ is the constant in Lemma \ref{s2:GN2}, then for any $(x,t),(y,s)\in Q'_{\phi}(x_0,h)$ we have
\begin{eqnarray*}
|u(x,t)-u(y,s)|\leq C^*h^{-\bz_0}\l(|x-y|+|t-s|\r)^{\alpha_0},
\end{eqnarray*}
where $\bz_0>0, \alpha_0\in(0,1)$ depend only on $n$; the constants $C^*,h_0$ depend only on $n,\Omega',T$.
\begin{proof}
In this proof we denote by $C,c$ constants depending only on $n,\Omega',T$. Their values may change from line to line whenever there is no confusion. Let $c\in(0,1)$ be a small constant to be chosen. Fix $h\leq c$ and $(x,t),(y,s)\in Q'_\phi(x_0,h)$ with $s\leq t$. Consider these two cases:\\

$\mathbf{Case\;1.}$ $(y,s)\in B(x,ct)\times (\f{7t}{8},t]$. If $c$ is small then by Lemma \ref{s2:GN1} we have,
\begin{equation*}
y\in\overline{B(x,|y-x|)}\subset\overline{S_\phi(x,\rz_0)},\quad\rz_0=C|x-y|.
\end{equation*}
If $c\le\f{1}{8C}$ is small, then $\rz_0\le \f{t}{8}$.

We have $(y,s)\in Q_\phi\l((x,t),\rz_0+t-s\r)$ with $\rz_0+t-s\leq\f{t}{4}$. Apply Theorem \ref{s2:H3} with $z_0\rightsquigarrow (x,t), R\rightsquigarrow t$ and $\rz\rightsquigarrow\rz_0+t-s$ and $u\rightsquigarrow u-u(x,0)$, we obtain
\begin{eqnarray}\label{s4:lem29}
|u(x,t)-u(y,s)|&\leq&osc_{Q_\phi((x,t),\rz_0+t-s)}u\nonumber\\
&\leq&C\l(\f{\rz_0+t-s}{t}\r)^{\alpha}\l[\sup_{Q_\phi((x,t),t)}|u-u(x,0)|+t^{\f{n}{2(n+1)}}\|f\|_{L^{n+1}(Q_\phi((x,t),t))}\r]\nonumber\\
&\leq&C\l(\f{\rz_0+t-s}{t}\r)^{\alpha}\l[\sup_{Q_\phi((x,t),t)}|u-u(x,0)|+t^2\|f\|_{\Lfz(Q_\phi((x,t),t))}\r],
\end{eqnarray}
where $\az=\az(n)$ and we use Lemma \ref{s2:GN3} in the last equality. Apply Lemma \ref{s4:lem1} and then for any $(x',t')\in Q_\phi((x,t),t)\subset Q'_\phi (x,h)\subset Q'_{\phi}(x_0,\tz_0^2h)$ ($\tz_0$ is the constant in Lemma \ref{s2:GN2}) we have
\begin{eqnarray}\label{s4:lem211}
|u(x',t')-u(x,0)|\leq Ch^{-\f{b}{1-b}}(\phi_x(x')+t')^{b}\leq C h^{-\f{b}{1-b}}t^{b},
\end{eqnarray}
where $b=b(n)$. Combining \eqref{s4:lem29}, \eqref{s4:lem211} and taking $\alpha_0:=\min\{\alpha,b\},\bz_0:=\f{b}{1-b}$, we obtain
\begin{eqnarray}\label{s4:lem212}
|u(x,t)-u(y,s)|\leq Ch^{-\f{b}{1-b}}\l(\f{\rz_0+t-s}{t}\r)^{\alpha}t^{b}\leq Ch^{-\bz_0}(|x-y|+t-s)^{\alpha_0}.
\end{eqnarray}

$\mathbf{Case\;2}$. $y\notin B(x,ct)$, or $\;t-s\geq\f{t}{8}$. Applying Lemma \ref{s4:lem1} we obtain
\begin{eqnarray*}
|u(x,t)-u(x,0)|&\leq&Ch^{-\f{b}{1-b}}(\phi_x(x)+t)^{b}\leq Ch^{-\f{b}{1-b}}t^{b}\leq Ch^{-\f{b}{1-b}}[|x-y|+|t-s|]^{b}.
\end{eqnarray*}
The estimate of $|u(y,s)-u(y,0)|$ is similar. Hence we obtain
\begin{eqnarray}\label{s4:lem213}
|u(x,t)-u(y,s)|&\leq&Ch^{-\bz_0}(|x-y+|t-s|)^{\alpha_0}.
\end{eqnarray}

Combining \eqref{s4:lem212} and \eqref{s4:lem213}, the proof of Lemma \ref{s4:lem2} is complete.
\end{proof}
\end{lem}

Since under the assumption $\mathbf{(H')}$ the operator $\mathcal{L}_{w}$ is uniformly parabolic in $\overline{B_1}\times[0,1]$, where $w$ is the solution of \eqref{w}, we have pointwise $C^{1+\az,\f{1+\az}{2}}$ regularity estimate for solutions $h$ of $\mathcal{L}_wh=0$ on the initial surface as follows.

\begin{lem}\label{s4:lem3}
Let $B_{\f{6}{5}}\subset\Omega\subset B_n$ be a normalized convex domain and $w$ is the convex solution of \eqref{w}. Assume $h:\overline{B_1}\times [0,1]\to\mathbb{R}$ is a continuous function satisfying $\|h\|_{\Lfz(B_1\times (0,1])}\le 1$ and
$$\left\{
\begin{array}{rcl}
\mathcal{L}_w h=0&&{\mathrm{in}\;B_1\times (0,1]},\\
h=\varphi&&{\mathrm{on}\;\overline{B_1}\times \{0\}}
\end{array}\right.$$
with
$$
|\varphi(x,0)-l_{\varphi,(0,0)}(x)|\leq |x|^{2}\quad\forall x\in\overline{B_1},\quad\quad\mathrm{and}\quad\quad |D\varphi(0,0)|\leq 1.
$$
Then for any $(x,t)\in B_1\times [0,1]$, we have
\begin{eqnarray*}
|h(x,t)-l_{\varphi,(0,0)}(x)|\leq K(n)(|x|^2+t).
\end{eqnarray*}
\end{lem}

Next we prove the regularity of solutions of \eqref{u} at the minimum point of $\phi$ on the initial surface, using the perturbation arguments as in the interior estimate.

\begin{thm}\label{s4:thm1}
Assume the hypothesis $\mathbf{(H)}$ holds with $\lambda=1-\tz,\Lambda=1+\tz$. Let $u$ be a continuous solution of \eqref{u} with $0<T_1<T<T_2$ and
$$\|u\|_{\Lfz(D)}+\|f\|_{\Lfz(D)}+\|\varphi(\cdot,0)\|_{C^{1,1}(\overline{\Omega})}\le 1.$$
Then for any $0<\alpha<1$ and $(x,t)\in D$ we have
\begin{eqnarray*}
|u(x,t)-l_{\varphi,(x_0,0)}(x)|\leq C(|x-x_0|^2+t)^{\f{1+\alpha}{2}},
\end{eqnarray*}
where $\tz\in(0,1),C>0$ depend only on $n,\alpha,T_1,T_2$ and $x_0$ is the minimum point of $\phi$.

\begin{proof}
After dividing $u$ by the constant $3+\tz^{-1}$ we may assume that
\begin{eqnarray}\label{s4:thm12}
2\|u\|_{\Lfz(D)}+\tz^{-1}\|f\|_{\Lfz(D)}+\|\varphi(\cdot,0)\|_{C^{1,1}(\overline{\Omega})}\leq 1.
\end{eqnarray}
Assume for simplicity that $x_0=0$. Denote $S_{\!\mu}(\phi):=S_{\!\phi}(0,\mu),\;Q_\mu'(\phi):=S_{\!\mu}(\phi)\times (0,\mu]$ and $l(x):=l_{\varphi,(0,0)}(x)$. Let $\alpha_1\in(0,1)$ be a constant to be chosen later. We will prove by introduction that there exist $\mu\in(0,1)$ depending only on $n,\alpha_1,\alpha,T_1,T_2$, a sequence of positive matrices $A_k$ with $\mathrm{det}\,A_k=1$ such that if we denote $\I_k(x,t):=(A_kx,t)$, then
\begin{equation*}
\|A_{k-1}A_k^{-1}\|\leq\f{1}{\sqrt{c_0}},\quad\|A_k\|\leq\sqrt{C_0(1+C_0\dz_0)(1+C_0\dz_1)\cdots(1+C_0\dz_{k-1})};\tag{a}
\end{equation*}
\begin{equation*}
\|A_k^{-1}\|\leq\f{1}{\sqrt{c_0(1-C_0\dz_0)(1-C_0\dz_1)\cdots(1-C_0\dz_{k-1})}};\tag{b}
\end{equation*}
\begin{equation*}
B_{(1-\dz_k)\sqrt{2}}(0)\times (0,1]\subset\mu^{-\f{k}{2}}\I_k Q'_{\mu^k}(\phi)\subset B_{(1+\dz_k)\sqrt{2}}(0)\times (0,1];\tag{c}
\end{equation*}
\begin{equation*}
\|u-l\|_{\Lfz(Q'_{\mu^k}(\phi))}\leq\mu^{\f{1+\alpha_1}{2}(k-1)};\tag{d}
\end{equation*}
\begin{equation*}
\f{|(u-l)(\mu^{\f{k}{2}}\I_k^{-1}(x,t))-(u-l)(\mu^{\f{k}{2}}\I_k^{-1}(y,s))|}{\mu^{\f{1+\alpha_1}{2}(k-1)}}
\leq 3 C^*\mu^{-\bz_0}(|x-y|+|t-s|^{\f{1}{2}})^{\alpha_0},\tag{e}
\end{equation*}
for any $(x,t),(y,s)\in\mu^{-\f{k}{2}}\I_k Q'_{\mu^k}(\phi)$, where $c_0,C_0$ are the constants in Lemma \ref{s2:GN5}, $C^*,\alpha_0$ and $\bz_0$ are given in Lemma \ref{s4:lem2}, the parabolic dilation is with respect to $(0,0)$ and is defined as in \eqref{d}, and\\

$A_0:=I,\quad \dz_0:=0,\quad \dz_1:=C_0(\mu^{\f{1}{2}}+\mu^{-1}\tz^{\f{1}{2}}),$

$\dz_k:=C_0(\dz_{k-1}\mu^{\f{1}{2}}+\mu^{-1}\tz^{\f{1}{2}}),\quad \dz_k<\dz_{k-1}\quad \mathrm{for}\; k\geq 2.$\\

The proof is similar to that of \cite[Theorem 4.5]{GN}. We just sketch it. The conclusion for $k=1$ easily follows from Lemma \ref{s2:GN5} (i) and Lemma \ref{s4:lem2}.
Suppose that we have (a)-(e) for $k=i\geq 1$, now we will verify (d) for $k=i+1$ and then construct $A_{i+1},\I_{i+1}$ and verify (a), (b), (c), (e) for $k=i+1$. 

+ Verifying (d) for $k=i+1$: First it is easily seen that
\begin{eqnarray*}
\|A_i\|, \|A_i^{-1}\|\leq C_1\exp\Big\{\f{2 C_0^2\sqrt{\tz}}{\mu}i\Big\},
\end{eqnarray*}
where $C_1$ is a constant depending only on $c_0,C_0$ and $\mu$. Let $0<\alpha_2<1$ be a constant to be chosen later, and take $\tz$ such that
$$\sqrt{\tz}\leq\alpha_2\f{\mu\ln\mu^{-1}}{4 C_0^2},$$
then
\begin{eqnarray}\label{s4:thm17}
\|A_i\|,\|A_i^{-1}\|\leq C_1\mu^{-\f{\alpha_2}{2}i}.
\end{eqnarray}
Let
\begin{eqnarray*}
\phi^*(y):&=&\f{1}{\mu^i}[\phi(\mu^{\f{i}{2}}A_i^{-1}y)-\phi(0)-\mu^i],\quad \mathrm{y}\in \Omega_i^*:=\mu^{-\f{i}{2}}A_iS_{\!\mu^i}(\phi),\\
v(y,s):&=&\f{1}{\mu^{\f{1+\alpha_1}{2}(i-1)}}(u-l)(\mu^{\f{i}{2}}A_i^{-1}y,\mu^i s),\quad \mathrm{(y,s)}\in\Omega_i^*\times (0,1].
\end{eqnarray*}
Then 
$$\mathcal{L}_{\phi^*} v(y,s)=f^*(y,s):=\f{\mu^i}{\mu^{\f{1+\alpha_1}{2}(i-1)}}f(\mu^{\f{i}{2}}A_i^{-1}y,\mu^i s)\quad \mathrm{in}\; \Omega_i^*\times (0,1].$$
Let $\epsilon>0$ be a constant to be chosen later. We apply Lemma \ref{s3:lem4.2} and we choose $\tz\leq\dz$, where $\dz>0$ is given by Lemma \ref{s3:lem4.2}. Let $h$ be a solution of
$$\left\{
\begin{array}{rcl}
\mathcal{L}_w h=0&&{\mathrm{in}\;B_1\times (0,1]},\\
h=v&&{\mathrm{on}\;\partial_p(B_1\times (0,1])},
\end{array}\right.$$
where $w$ is the convex solution of the equation
$$\left\{
\begin{array}{rcl}
\mathrm{det} D^2 w=1&&{\mathrm{in}\;\Omega_i^*},\\
w=0&&{\mathrm{on}\;\partial\Omega_i^*},
\end{array}\right.$$
then we have
\begin{eqnarray}\label{s4:thm18}
\|v-h\|_{\Lfz\l(B_1\times (0,1]\r)}\leq\epsilon.
\end{eqnarray}
Since
\begin{eqnarray*}
h(y,0)=\varphi^*(y):=\f{1}{\mu^{\f{1+\alpha_1}{2}(i-1)}}[\varphi(\mu^{\f{i}{2}}A_i^{-1}y,0)-l_{\varphi,(0,0)}(\mu^{\f{i}{2}}A_i^{-1}y)],
\end{eqnarray*}
and note that 
\begin{eqnarray}\label{varphi^*}
\|\varphi^*\|_{C^{1,1}(\overline{\Omega_i^*})}\leq\f{1}{\mu^{\f{1+\alpha_1}{2}(i-1)}}\|\mu^{\f{i}{2}}A_i^{-1}\|^{2}
\leq\f{(\mu^{\f{i}{2}}C_1\mu^{-\f{\alpha_2}{2}i})^{2}}{\mu^{\f{1+\alpha_1}{2}(i-1)}}\leq \mu^{\f{1+\az_1}{2}},
\end{eqnarray}
where we choose $\alpha_2$ and $\alpha_1$ such that $2(1-\alpha_2)>1+\alpha_1$. Then by Lemma \ref{s4:lem3},
\begin{eqnarray*}
|h(y,s)|\leq 2K(|y|^2+s)\quad\quad\forall (y,s)\in B_1\times (0,1],
\end{eqnarray*}
where $K=K(n)$. Lemma \ref{s2:GN5} (i) implies that
\begin{eqnarray*}
|h(y,s)|\leq 6KC_0^2\tz_0^2\mu\quad\quad\forall (y,s)\in Q'_{\tz_0^2\mu}(\phi^*),
\end{eqnarray*}
where we recall $Q'_{\tz_0^2\mu}(\phi^*):=S_{\!\phi^*}(0,\tz_0^2\mu)\times(0,\tz_0^2\mu]$. Choose $\epsilon:=KC_0^2\tz_0^2\mu$, and combine \eqref{s4:thm18} with the last estimate we get
\begin{eqnarray}\label{s4:thm110}
|v(y,s)\leq 8KC_0^2\tz_0^2\mu\leq\mu^{\f{1+\alpha_1}{2}}\quad\quad\forall (y,s)\in Q'_{\tz_0^2\mu}(\phi^*).
\end{eqnarray}
Back to $u$ we find that $(d)$ for $k=i+1$ holds.

+ Constructing $A_{i+1}$: We apply Lemma \ref{s2:GN5} (ii) to $\Omega^*_i$ and obtain that (a)-(c) for $k=i+1$ hold. Combining \eqref{varphi^*} and \eqref{s4:thm110}, we can apply Lemma \ref{s4:lem2} and find that $(e)$ for $k=i+1$ holds.\\

By \eqref{s4:thm17}, we have for each $k\geq 1$,
$$B(0, c\mu^{\f{1+\alpha_2}{2}k})\times (0,c^2\mu^{(1+\alpha_2)k}]\subset Q'_{\mu^k}(\phi)$$
for some constant $c$ depending only on $C_1$. Given $\alpha\in(0,1)$, we choose $\alpha_2$ small such that
$$2(1-\alpha_2)>1+(1+\alpha_2)\alpha+\alpha_2.$$ Take $\alpha_1:=(1+\alpha_2)\alpha+\alpha_2$, which satisfies
$$2(1-\alpha_2)>1+\alpha_1,$$
as required in the proof, then we have $\f{1+\alpha_1}{1+\alpha_2}=1+\alpha$, it follows that for each $k_0\geq 1$ and any $(x,t)\in B(0,c\mu^{\f{1+\alpha_2}{2}k_0})
\times(0,c^2\mu^{(1+\alpha_2)k_0}]$,
\begin{eqnarray*}
|u(x,t)-l(x)|\leq\mu^{\f{1+\alpha_1}{2}(k_0-1)}=c^{-(1+\alpha)}\mu^{-\f{1+\alpha_1}{2}}[c\mu^{\f{1+\alpha_2}{2}k_0}]^{1+\alpha}.
\end{eqnarray*}
\end{proof}
\end{thm}

\begin{rem}\label{s4:rem1}
In the above theorem it can also be concluded from (d) that
\begin{eqnarray*}
\max_{\overline{Q'_r(\phi)}}|u(x,t)-l(x)|\leq C r^{\f{1+\alpha}{2}},\quad\quad\forall r>0.
\end{eqnarray*}
where $C$ depends only on $n,\alpha,T_1,T_2$.
\end{rem}

Using Theorem \ref{s4:thm1} we can obtain pointwise regularity of solutions of \eqref{u} near the initial surface when the density of the Monge-Amp$\grave{e}$re measure is continuous.

\begin{thm}\label{s4:thm2}
Assume $\mathbf{(H)}$ holds and $g:=\mathrm{det} D^2\phi\in C(\Omega)$. Assume that $\Omega'\Subset\Omega,\;0<\alpha<1$, and $u$ is a solution of \eqref{u} with
$$\|u\|_{\Lfz(D)}+\|f\|_{\Lfz(D)}+\|\varphi(\cdot,0)\|_{C^{1,1}(\overline{\Omega})}\le 1,$$
then for any $x_0\in\Omega'$ and $(x,t)\in D$ we have
\begin{eqnarray*}
|u(x,t)-l_{\varphi,(x_0,0)}(x)|&\leq& C(|x-x_0|^2+t)^{\f{1+\alpha}{2}},
\end{eqnarray*}
where $C>0$ depends only on $n,\lambda,\Lambda,\alpha,\Omega',T$ and the modulus of continuity of $g$.
\end{thm}

\begin{rem}\label{s4:rem2}
Under the assumptions of Theorem \ref{s4:thm2}, we can also obtain from Remark \ref{s4:rem1} that 
\begin{eqnarray*}
\max_{\overline{Q'_{\phi}(x_0,r)}}|u(x,t)-l_{\varphi,(x_0,0)}(x)|\leq C r^{\f{1+\alpha}{2}},\quad\quad\forall x_0\in\Omega',\forall r>0,
\end{eqnarray*}
where $C>0$ depends only on $n,\lambda,\Lambda,\alpha,\Omega'$, $T$ and the modulus of continuity of $g$.
\end{rem}

Next we prove $C^{1+\az,\f{1+\az}{2}}$ estimates for solutions of \eqref{u} near the initial surface.

\begin{thm}\label{s4:thm3}
Assume the hypotheses of Theorem 4.2 hold. Then $u\in C^{1+\alpha,\f{1+\alpha}{2}}(\overline{\Omega'}\times [0,T])$ for any $\alpha\in(0,1)$ and
\begin{eqnarray*}
\|u\|_{C^{1+\alpha,\f{1+\alpha}{2}}(\overline{\Omega'}\times [0,T])}\leq C,
\end{eqnarray*}
where $C>0$ depend only on $n,\lambda,\Lambda,\alpha,\Omega',T$ and the modulus of continuity of $g$. 

\begin{proof}
Fix any $\az'\in(0,1)$. In this proof we denote by $C,c$ constants depending only on $n,\lambda,\Lambda,\az',\Omega',T$ and the modulus of continuity of $g$.

Fix $x_0\in\Omega'$ and $t_0$ small. Let $Tx=A(x-x_0)+y_0$ be an affine transformation such that $B_1\subset T S_{\!\phi}(x_0,t_0)\subset B_n$. Denote $T_{\!p}(x,t):=(Tx,t_0^{-1}t)$. By interior $C^{1,\gamma}$ estimates for $\phi$,\begin{equation}\label{s4:thm31'}
B(x_0,c\,t_0^{\f{1}{1+\gamma}})\subset S_{\!\phi}(x_0,t_0)
\end{equation}
for some $\gz=\gz(n,\lz,\Lz)\in(0,1)$, which gives that
\begin{equation}\label{s4:thm31}
\|A\|\leq C\,t_0^{-\f{1}{1+\gamma}}.
\end{equation}

We point out that we can choose $\gz$ close to $1$. Indeed, since $g\in C(\Omega)$, by the interior $W^{2,p}$ estimates for solutions of Monge-Amp$\grave{e}$re equations in \cite{C}, we have $\phi\in W^{2,p}_{\mathrm{loc}}$ for any $p<\infty$. Then the imbedding theorem \cite[Theorem 7.26]{GT} implies that $\phi\in C^{1,\gamma}_{\mathrm{loc}}$ for any $p>n$ and $\gamma<1-n/p$. Note that we may choose $p$ large such that $\gamma$ close to $1$.

Define
$$\phi^*(y)=\f{1}{t_0}[\phi(T^{-1}y)-l_{\phi,x_0}(T^{-1}y)-t_0],$$
$$u^*(y,s)=\f{1}{t_0}[u(T^{-1}y,t_0 s)-l_{\varphi,(x_0,0)}(T^{-1}y)].$$
Then 
\begin{eqnarray*}
&&B_1\times (0,1]\subset T_{\!p}(Q'_{\phi}(x_0,t_0))=Q'_{\phi^*}(y_0,1)\subset B_n\times (0,1],\\
&&\lambda'\leq\mathrm{det}D^2\phi^*(y)\leq\Lambda'\quad\mathrm{in}\;S_{\!\phi^*}(y_0,1);\quad\quad\phi^*=0\quad\mathrm{on}\;\partial S_{\!\phi^*}(y_0,1),\\
&&\mathcal{L}_{\phi^*} u^*(y,s)=f^*(y,s):=\f{1}{t_0^n(\mathrm{det}\,A)^2}f(T^{-1}y,t_0 s)\quad\mathrm{in}\;Q'_{\phi^*}(y_0,1),
\end{eqnarray*}
where $\lambda',\Lambda'>0$ depend only on $n,\lambda,\Lambda$.

Denote $D^*:=Q'_{\phi^*}(y_0,1)$. We apply the interior estimate Theorem $1$ and find that 
\begin{eqnarray}\label{s4:thm32}
\|u^*\|_{C^{1+\alpha',\f{1+\alpha'}{2}}\l(\overline{S_{\phi^*}(y_0,\f{1}{2})}\times [\f{1}{2},1]\r)}\leq C[\|u^*\|_{\Lfz(D^*)}+\|f^*\|_{\Lfz(D^*)}].
\end{eqnarray}

We take $\az>\az'$ and then choose $\gz$ close to $1$ such that
\begin{eqnarray}\label{s4:thm3*}
\f{1+\alpha}{2}>\f{1+\az'}{1+\gamma}.
\end{eqnarray}
From Remark \ref{s4:rem2} we obtain that 
\begin{eqnarray}\label{s4:thm33}
\|u^*\|_{\Lfz(D^*)}\leq C\f{1}{t_0}t_0^{\f{1+\alpha}{2}}.
\end{eqnarray}
For any $z_i=(x_i,t_i)\in S_{\!\phi}\l(x_0,\f{t_0}{2}\r)\times [\f{t_0}{2},t_0], i=1,2$, the estimates \eqref{s4:thm31}-\eqref{s4:thm33} imply that
\begin{eqnarray}\label{s4:thm34}
|D u(z_1)&-&D u(z_2)|=t_0|A^t(D u^*(T x_1,t_0^{-1}t_1)-D u^*(T x_2,t_0^{-1}t_2))|\nonumber\\
&\leq&Ct_0^{\f{1+\alpha}{2}-\f{1+\az'}{1+\gamma}}(|x_1-x_2|+|t_1-t_2|^{\f{1}{2}})^{\alpha'}\nonumber\\
&\leq&C(|x_1-x_2|+|t_1-t_2|^{\f{1}{2}})^{\alpha'}.
\end{eqnarray}
Similarly, we have
\begin{eqnarray}\label{s4:thm34'}
|u(x_1,t_1)-u(x_1,t_2)|&=&t_0|u^*(T x_1,t_0^{-1}t_1)-u^*(T x_1,t_0^{-1}t_2)|\nonumber\\
&\leq&Ct_0^{\f{\alpha-\alpha'}{2}}|t_1-t_2|^{\f{1+\alpha'}{2}}\le C|t_1-t_2|^{\f{1+\alpha'}{2}},
\end{eqnarray}
and
\begin{eqnarray}\label{s4:thm34''}
|D u(x_1,t_1)|\leq t_0|A^t D u^*(T x_1,t_0^{-1}t_1)|+|D\varphi(x_0,0)|\leq Ct_0^{\f{1+\alpha}{2}-\f{1}{1+\gamma}}\le C.
\end{eqnarray}
The estimates \eqref{s4:thm34} and \eqref{s4:thm34'} imply that for any $(x,t)\in S_{\!\phi}\l(x_0,\f{t_0}{2}\r)\times [\f{t_0}{2},t_0]$,
\begin{eqnarray}\label{s4:thm35}
|u(x,t)-l_{u,(x_0,t_0)}(x)|\leq C\l(t_0^{\f{\alpha-\alpha'}{2}}|t-t_0|^{\f{1+\alpha'}{2}}+|x-x_0|^{1+\alpha'}\r).
\end{eqnarray}

On the other hand, from Theorem \ref{s4:thm2}, for any $(x,t)\in D$ we have
\begin{eqnarray}\label{s4:thm35'}
|u(x,t)-l_{\varphi,(x_0,0)}(x)|
\leq C(|x-x_0|^2+t)^{\f{1+\alpha}{2}}.
\end{eqnarray} 

By \cite[Lemma 4.8]{GN}, there exists a function $\psi\in C^\infty_0(\mathbb{R}^n)$ with support in the unit ball such that $\psi_\epsilon\ast l=l$ for each $\epsilon>0$ and every linear function $l$. As usual, $\psi_\epsilon(x):=\epsilon^{-n}\psi(\f{x}{\epsilon})$. Similar to \cite[Page 2068-2069]{GN} we can write
\begin{eqnarray*}
D_j u(x_0,t_0)-D_j\varphi(x_0,0)=[u(\cdot,t_0)-l_{\varphi,(x_0,0)}]\ast D_j\psi_\epsilon (x)
-[u(\cdot,t_0)-l_{u,(x_0,t_0)}]\ast D_j\psi_\epsilon (x).
\end{eqnarray*}
Set $x:=x_0$, and $\epsilon:=c\l(\f{t_0}{2}\r)^{\f{1}{1+\gamma}}$, then by \eqref{s4:thm35} we have
\begin{eqnarray*}
|[u(\cdot,t_0)-l_{u,(x_0,t_0)}]\ast D_j\psi_\epsilon (x_0)|
&=&\epsilon^{-n-1}\l|\int_{B_\epsilon(x_0)}[u(y,t_0)-l_{u,(x_0,t_0)}(y)]D_j\psi\l(\f{x_0-y}{\epsilon}\r)dy\r|\\
&\leq&C\epsilon^{-n-1}\|D_j\psi\|_{\infty}\int_{|y-x_0|\leq \epsilon}|y-x_0|^{1+\alpha'}dy\\
&\leq&Ct_0^{\f{\alpha'}{1+\gamma}},
\end{eqnarray*}
also, by \eqref{s4:thm35'},
\begin{eqnarray*}
|[u(\cdot,t_0)-l_{\varphi,(x_0,0)}]\ast D_j\psi_\epsilon (x_0)|
&=&\epsilon^{-n-1}\l|\int_{B_\epsilon(x_0)}[u(y,t_0)-l_{\varphi,(x_0,0)}(y)]D_j\psi\l(\f{x_0-y}{\epsilon}\r)dy\r|\\
&\leq&C\epsilon^{-n-1}\|D_j\psi\|_{\infty}\int_{B_\epsilon(x_0)}(|y-x_0|^2+t_0)^{\f{1+\alpha}{2}}dy\\
&\leq&C\l(t_0^{\f{1+\alpha}{2}-\f{1}{1+\gamma}}+t_0^{\f{\alpha}{1+\gamma}}\r),
\end{eqnarray*}
combining these inequalities we obtain
\begin{eqnarray}\label{s4:thm37}
|D u(x_0,t_0)-D\varphi(x_0,0)|&\leq&C t_0^{\f{\alpha'}{1+\gamma}}.
\end{eqnarray}
Combining \eqref{s4:thm34}, \eqref{s4:thm34'}, \eqref{s4:thm34''}, \eqref{s4:thm35'} and \eqref{s4:thm37}, if we denote 
$$Q_{x,t}:=B(x,c\,(t/2)^{\f{1}{1+\gz}})\times [t/2,t],$$
then the estimate
\begin{eqnarray*}
\|u\|_{C^{1+\alpha',\f{1+\alpha'}{2}}(\overline{Q_{x_0,t_0}})}+\f{|u(x_0,t_0)-\varphi(x_0,0)|}{t^{\f{1+\alpha}{2}}}+\f{|D u(x_0,t_0)-D\varphi(x_0,0)|}{t_0^{\f{\alpha'}{1+\gamma}}}\quad\leq C.
\end{eqnarray*}
holds for any $(x_0,t_0)$ such that $x_0\in\Omega'$ and $t_0$ small. This  implies that
\begin{eqnarray}\label{s4:thm310}
\|u\|_{C^{1+\alpha',\f{1+\alpha'}{2}}(\overline{\Omega'}\times[0,c])}\leq C.
\end{eqnarray}
Recall \eqref{s4:thm3*} for the definition of $\alpha'$. Combining \eqref{s4:thm310} and the interior estimate Theorem 1, the proof of Theorem \ref{s4:thm3} is complete.
\end{proof}
\end{thm}

\section{Regularity near the side of $\Omega$}\label{s5}

In this section, we fix constants $0<\lz\le\Lz<\infty,\rz>0$ and refer to all positive constants depending only $n,\lz,\Lz$ and $\rz$ as universal constants.

We assume the following condition on $\Omega$ and $\phi$: Assume $\Omega\subset\mathbb{R}^n$ is a bounded convex set with
\begin{equation}\label{s5:lem11}
B_\rz(\rz e_n)\subset\Omega\subset \{x_n\geq 0\}\cap B_{\f{1}{\rz}}
\end{equation}
and
\begin{equation}\label{s5:lem12}
\mathrm{\Omega\;contains\;an\;interior\;ball\;of\;radius\;\rz\;
tangent\;to\;\partial\Omega\;at\;each\;point\;on\;\partial\Omega\cap\;B_\rz}.
\end{equation}
Let $\phi:\overline{\Omega}\rightarrow\mathbb{R}$, $\phi\in C^{0,1}(\overline{\Omega})\cap C^2(\Omega)$ be a convex function satisfying
\begin{equation}\label{s5:lem13}
\mathrm{det} D^2\phi=g,\quad\quad 0<\lambda\leq g\leq\Lambda\quad\mathrm{in}\;\Omega.
\end{equation}
Assume further that on $\partial\Omega\cap B_\rz$, $\phi$ separates quadratically from its tangent planes on $\partial\Omega$, namely, for any $x_0\in\partial\Omega\cap B_\rz$ we have
\begin{equation}\label{s5:lem14}
\rz|x-x_0|^2\leq\phi(x)-\phi(x_0)-\nabla\phi(x_0)\cdot (x-x_0)\leq\rz^{-1}|x-x_0|^2,\quad\quad\forall x\in\partial\Omega.
\end{equation}

Let $0<c_*<T$. We establish regularity of $u$ near $\partial\Omega$ where $u:\l(B_\rz\cap\overline{\Omega}\r)\times [\f{c_*}{2},T]\rightarrow\mathbb{R}$ is a continuous solution of
\begin{equation}\label{us}
\left\{
\begin{array}{rcl}
\L u=f&&{\mathrm{in}\;(B_\rz\cap\Omega)\times (\f{c_*}{2},T]},\\
u=0&&{\mathrm{on}\;(B_\rz\cap\partial\Omega)\times [\f{c_*}{2},T]}.
\end{array}\right.
\end{equation}

The arguments in \cite{S1} easily apply to the parabolic case as long as we construct corresponding supersolution and subsolution and employ the weak Harnack inequality \cite{HQ} for nonnegative supersolutions of $\mathcal{L}_\phi u=0$. Hence we just sketch the proof of the results which are straightforward modifications of the elliptic case \cite{S1}, and give more detail whenever necessary. 

We use frequently the two localization theorems below concerning geometry of boundary sections and maximal interior sections of $\phi$ respectively.

\begin{thm}(Localization Theorem \cite{S1,S2})\label{s5:thm1}
Assume $\Omega$ satisfies \eqref{s5:lem11} and $\phi$ satisfies \eqref{s5:lem13}, and
$$\phi(0)=\nabla\phi(0)=0,\quad\quad
\rz|x|^2\leq\phi(x)\leq\rz^{-1}|x|^2\quad\mathrm{on}\;\partial\Omega\cap\{x_n\leq\rz\}.$$
Then there exists a universal constant $k$ such that for each $h\leq k$, there is an ellipsoid $E_h$ of volume $|B_1|h^{n/2}$ satisfying
$$k E_h\cap\overline{\Omega}\subset S_{\!\phi}(0,h)\subset k^{-1}E_h.$$
Moreover, the ellipsoid $E_h$ is obtained from the ball of radius $h^{\f{1}{2}}$ by a linear transformation $A_h^{-1}$ (sliding along the $x_n=0$ plane)
\begin{eqnarray*}
\mathrm{det}\,A_h=1,\quad A_h x=x-\tau_h x_n,\quad \tau_h\cdot e_n=0,\\
h^{-\f{1}{2}}A_h E_h=B_1,\quad |\tau_h|\leq k^{-1}|\log h|.
\end{eqnarray*}
\end{thm}

\begin{prop}(See \cite[Proposition 3.2]{S1}.)\label{s5:prop1}
Let $\phi$ and $\Omega$ satisfy the hypotheses of the Localization Theorem \ref{s5:thm1}. Assume that for some $y\in\Omega$ the section $S_{\!\phi}(y,h)\subset\Omega$ is tangent to $\partial\Omega$ at $0$ for some $h\leq c$ with $c$ universal. Then 
$$\nabla\phi(y)=a e_n\quad\quad\mathrm{for\;some\;}a\in [k_0 h^{\f{1}{2}},k_0^{-1}h^{\f{1}{2}}],$$
$$k_0 E_h\subset S_{\!\phi}(y,h)-y\subset k_0^{-1}E_h,\quad\quad k_0 h^{\f{1}{2}}\leq\mathrm{dist}(y,\partial\Omega)\leq k_0^{-1}h^{\f{1}{2}},$$
where $E_h$ is the ellipsoid defined in the Localization Theorem \ref{s5:thm1} and $k_0$ is a universal constant.
\end{prop}

Under the assumptions of Theorem \ref{s5:thm1} we have
\begin{equation}\label{s5:lem14''}
\overline{\Omega}\cap B_{c_1 h^{\f{1}{2}}/|\log h|}^+\subset S_{\!\phi}(0,h)\subset B_{C_1 h^{\f{1}{2}}|\log h|},
\end{equation}
and the estimate
\begin{equation}\label{s5:lem14'}
c_1|x|^2|\log|x||^{-2}\leq\phi(x)\leq C_1|x|^2|\log|x||^{2},
\end{equation}
holds in a neighborhood of $0$, where $c_1,C_1$ are universal constants (see Equation (4.3) in \cite{S1}).

In the lemma below, we construct a supersolution and estimate $\f{u(\cdot,t)}{d_{\partial\Omega}}$ near $\partial\Omega$.

\begin{lem}\label{s5:lem1}
Assume \eqref{s5:lem11}-\eqref{s5:lem14} hold. Let $u$ be a continuous solution of \eqref{us} with
$$\|u\|_{\Lfz\l((B_\rz\cap\Omega)\times (\f{c_*}{2},T]\r)}+\|f/\mathrm{tr}\,\Phi\|_{\Lfz\l((B_\rz\cap\Omega)\times (\f{c_*}{2},T]\r)}\le 1.$$
Then 
\begin{eqnarray*}
|u(x,t)|\leq Cd_{\partial\Omega}(x)\quad\quad\mathrm{in}\;(B_{c_0}\cap\Omega)\times [c_*,T],
\end{eqnarray*}
where $c_0$ is a universal constant and $C$ depends only on $n,\lambda,\Lambda,\rz,c_*$.

\begin{proof}
Let $\tilde{\delta}>0$ be chosen later, after multiplying $u$ by $\tilde{\delta}$ we may assume
\begin{equation*}\label{s5:lem15}
\|u\|_{\Lfz\l((B_\rz\cap\Omega)\times (\f{c_*}{2},T]\r)}+\|f/\mathrm{tr}\,\Phi\|_{\Lfz\l((B_\rz\cap\Omega)\times (\f{c_*}{2},T]\r)}\leq\tilde{\delta}.
\end{equation*}
Denote
$$q(x):=\f{1}{2}\l(\tilde{\delta}|x'|^2+\f{\Lambda^n}{(\lambda\tilde{\delta})^{n-1}}x_n^2\r)$$
and
$$\bar{w}_{t_0}(x,t):=M x_n-n(t-t_0)+\phi(x)-4 q(x)$$
for some large constant $M$ depending only on $n,\lambda,\Lambda,\rz,\tilde{\dz}$. Then it is straightforward to prove that 
$$|u|\leq\bar{w}_{t_0}\quad\quad\mathrm{in}\;(B_{c_0}\cap\Omega)\times [t_0-c_*/2,t_0]$$ 
if $\tilde{\delta}$ is small. Hence,
\begin{eqnarray*}
|u(0,x_n,t_0)|\leq M x_n+C_1 x_n^2|\log x_n|^{2}\leq C x_n,\quad\quad\forall x_n\in [0,c_0],
\end{eqnarray*}
where $C$ is a constant depending only on $n,\lambda,\Lambda,\rz,c_*$ and $C_1$ is as in \eqref{s5:lem14'}. The conclusion follows if we replace $0$ with each point $x_0\in B_{\f{\rz}{2}}\cap\partial\Omega$ and modify the construction of the corresponding supersolution.
\end{proof}
\end{lem}

Assume \eqref{s5:lem11}-\eqref{s5:lem14} hold with
$$\phi(0)=\nabla\phi(0)=0.$$ 
Let $E_h,A_h$ be as in the Localization Theorem \ref{s5:thm1}. Denote 
\begin{eqnarray*}
\phi_h(x):=\f{\phi(h^{\f{1}{2}}A_h^{-1}x)}{h}\quad\quad\mathrm{and}\quad\quad\Omega_h:=h^{-\f{1}{2}}A_h\Omega.
\end{eqnarray*}
Then
$$\mathrm{det} D^2\phi_h=g_h,\quad\quad 0<\lambda\leq g_h(x):=g(h^{\f{1}{2}}A_h^{-1}x)\leq\Lambda\quad\mathrm{in}\;\Omega_h,$$
$$h^{-\f{1}{2}}A_h S_{\!\phi}(0,h)=S_{\!\phi_h}(0,1):=\{x\in\overline{\Omega_h}:\phi_h(x)<1\},$$
$$\overline{\Omega_h}\cap B_k\subset S_{\!\phi_h}(0,1)\subset B_{k^{-1}}^+.$$

We introduce the class $\mathcal{D}_\sigma$ which consists of the triples $(\phi,\Omega,U)$ satisfying $(i)$-$(v)$ (See \cite{S1,LN}):

\noindent(i) $0\in\partial\Omega,\,U\subset\Omega\subset\mathbb{R}^n$ are bounded convex domains such that
$$B_k^+\cap\overline{\Omega}\subset\overline{U}\subset B_{k^{-1}}^+.$$
\noindent(ii) $\phi:\overline{\Omega}\rightarrow\mathbb{R}^+$ is convex satisfying $\phi=1\;\mathrm{on}\;\partial U\cap\Omega$ and
$$\phi(0)=0,\quad\nabla\phi(0)=0,\quad\lambda\leq\mathrm{det} D^2\phi\leq\Lambda\;\mathrm{in}\;\Omega,$$
$$\partial\Omega\cap\{\phi<1\}=\partial U\cap\{\phi<1\}.$$
\noindent(iii)
$$\f{\rz}{4}|x-x_0|^2\leq\phi(x)-\phi(x_0)-\nabla\phi(x_0)\cdot (x-x_0)\leq \f{4}{\rz}|x-x_0|^2,\quad\quad\forall x,x_0\in\partial\Omega\cap B_{\f{2}{k}}.$$
\noindent(iv)
$$\partial\Omega\cap\{\phi<1\}\subset G\subset\{x_n\leq\sigma\},$$
where $G\subset B_{2/k}$ is a graph in the $e_n$ direction and its $C^{1,1}$ norm is bounded by $\sigma$.

\noindent(v) $\phi$ satisfies in $U$ the hypotheses of the Localization Theorem in \ref{s5:thm1} at all points on $\partial U\cap B_{c_0}$ and if $r\leq c_0$, then
$$|\nabla\phi|\leq C_0 r|\log r|^2\quad\mathrm{in}\;\overline{\Omega}\cap B_r.$$
The constants $k,c_0,C_0$ above are universal.

It follows from \cite[Lemma 4.2]{S1} that if $h\leq h_0$, then $(\phi_h,\Omega_h,S_{\!\phi_h}(0,1))\in\mathcal{D}_\sigma$ with $\sigma=C h^{1/2}$, where $h_0,C>0$ are universal constants.

To iterate the estimate of $\f{u(\cdot,t)}{d_{\partial\Omega}}$ we need the lemma below deduced from the weak Harnack inequality \cite[Lemma 3.1, Theorem 3.1, Theorem 4.1]{HQ}.

\begin{lem}\label{s5:lem2}
Let $(\phi,\Omega,U)\in\mathcal{D}_\dz$ with $\dz$ small, universal. Let $u$ be a nonnegative classical solution of $\L u\leq 0$ in $Q=U\times(-1,0]$. Denote $z_0=(2\dz e_n,t_0)$ with $-\f{1}{2}<t_0<t_0+R\le 0$ and $R>0$ small.
\begin{enumerate}
\item[\rm(i)] If
\begin{equation*}
|\{z\in Q^-:u(z)\geq 1\}|\geq \mu|Q^-|
\end{equation*}
for some $0<\mu<1$, then
$$\inf_{Q^+}u\geq c,$$
where $Q^+=S_{\!\phi}(2\dz e_n,R)\times(t_0+R,t_0+2R]$, $Q^-=S_{\!\phi}(2\dz e_n,R)\times(t_0-R,t_0]$, and $c>0$ depends only on $n,\lambda,\Lambda,\rz$ and $\mu$.
\item[\rm(ii)] Suppose that
$$\inf_{S_{\!\phi}(2\dz e_n,R)\times(t_0,t_0+R]}u\leq 1,$$
then
$$\mathcal{M}(\{z\in Q_\phi(z_0,R): u(z)<M_2\})>\epsilon_0\mathcal{M}(Q_\phi(z_0,R)),$$
where $\mathcal{M}$ is the parabolic Monge-Amp$\grave{e}$re measure generated by $\phi(x)-t$, i.e., $d\mathcal{M}=\mathrm{det}D^2\phi dxdt$, and $M_2>0,\epsilon_0\in(0,1)$ are universal.
\end{enumerate}
\begin{proof}
We first prove $(i)$. Assume in contradiction that
$$\inf_{Q^+}u<c.$$
Then by \cite[Theorem 4.1]{HQ}, for any $i\geq 1$ we have
\begin{equation*}
|\{z\in Q^-:u(z)>c\tilde{K}M^i\}|\leq C\gamma^i|Q^-|,
\end{equation*}
where $\tilde{K},M,C>0,0<\gamma<1$ are universal constants. We reach a contradiction if $c$ is small.

For $(ii)$, by \cite[Theorem 3.1]{HQ}, we have
$$\mathcal{M}(\{z\in Q_\phi(z_0,R/2): u(z)<M_1\})>\epsilon_0\mathcal{M}(Q_\phi(z_0,R/2)),$$
which implies
$$\inf_{Q_\phi(z_0,R/2)}u\leq M_1.$$
Applying \cite[Lemma 3.1]{HQ} to $\f{u}{M_1}$, we obtain
$$\mathcal{M}(\{z\in Q_\phi(z_0,R): u(z)<M_0 M_1\})>\epsilon_0\mathcal{M}(Q_\phi(z_0,R)),$$
where $M_0,M_1$ are universal. Choose $M_2:=M_0 M_1$ and the conclusion follows.
\end{proof}
\end{lem}

Using the above lemma we can construct a subsolution and iterate the estimate of $\f{u(\cdot,t)}{d_{\partial\Omega}}$ in Lemma \ref{s5:lem1}.

\begin{lem}\label{s5:lem3}
Let $(\phi,\Omega,U)\in\mathcal{D}_\dz$ with $\dz$ small, universal. Denote $Q:=U\times (-1,0]$. Assume $u:\overline{Q}\to\mathbb{R}$ is a continuous solution of
$$\L u=f,\quad\quad a\,d_G\leq u\leq b\,d_G\quad\mathrm{in}\;Q,$$
where $a,b\in [-1,1]$ and
$$\max\{\dz,\|f/\mathrm{tr}\,\Phi\|_{\Lfz(Q)}\}\leq c_1(b-a).$$
Then
$$a'\,d_G\leq u\leq b'\,d_G\quad\mathrm{in}\;Q_\tz(\phi)=S_{\!\phi}(0,\tz)\times (-\tz,0],$$
where
$$a\leq a'\leq b'\leq b,\quad\quad b'-a'\leq\eta(b-a).$$
Here the constants $c_1,\eta\in(0,1),\tz$ are universal.

\begin{proof}
Set
$$u_1:=\f{u-a\,d_G}{b-a},\quad\quad u_2:=\f{b\,d_G-u}{b-a},$$
then
$$u_1,u_2\geq 0\quad\mathrm{in}\;Q,\quad\quad u_1+u_2=d_G.$$
If $R$ is small, then
$$d_G(y)\geq\dz,\quad\quad\forall y\in S_{\!\phi}(2\dz e_n,R).$$
Let $m\geq 1$ be such that
\begin{equation}\label{s5:lem3*}
m R<\dz^{\f{1}{n-1}}\leq(m+1)R.
\end{equation}
Denote 
$$Q^{(i)}:=S_{\!\phi}(2\dz e_n,R)\times (-(i+1)R,-iR], \quad i=0,1,\dots,m+2,$$ 
and assume without loss of generality that
\begin{eqnarray}\label{s5:lem32}
|\{u_1\geq\dz/2\}\cap Q^{(m+2)}|\geq\f{1}{2}|Q^{(m+2)}|.
\end{eqnarray}

Let $c_1>0$ be chosen later, and set $$\tilde{u}_1:=u_1+c_1(k^{-2}-|x|^2).$$ Then from $(\mathrm{i})$ in the definition of $\mathcal{D}_\dz$ we have $\tilde{u}_1\geq u_1\geq 0$. And it follows from $(\mathrm{iv})$ in the definition of $\mathcal{D}_\dz$ that $D^2 d_G\leq\dz I$, which gives
$$\L\tilde{u}_1=\f{\L u_1-a\,\mathrm{tr}(\Phi D^2 d_G)}{b-a}-2c_1\mathrm{tr}\,\Phi\leq 0,$$
this together \eqref{s5:lem32} and Lemma \ref{s5:lem2} $(i)$ implies that
\begin{equation*}
\inf_{Q^{(m)}}\tilde{u}_1\geq \widetilde{c_2}
\end{equation*}
for some universal constant $\widetilde{c_2}>0$.

Apply Lemma \ref{s5:lem2} $(ii)$ to $\f{M_2}{\widetilde{c_2}}\tilde{u}_1$ and we obtain
$$\inf_{Q^{m-1}}\f{M_2}{\widetilde{c_2}}\tilde{u}_1\geq 1.$$
Continue applying Lemma \ref{s5:lem2} $(ii)$ and noting \eqref{s5:lem3*}, we obtain
\begin{equation}\label{s5:lem33}
\tilde{u}_1\geq 2c_2:=\f{\widetilde{c_2}}{M_2^{\dz^{\f{1}{n-1}}/R}}\quad\quad\mathrm{in}\;S_{\!\phi}(2\dz e_n,R)\times (-\dz^{\f{1}{n-1}},0].
\end{equation}
The estimate \eqref{s5:lem33} holds for any $y\in F_\dz:=\{x_n=2\dz,\;|x'|\leq\dz^{\f{1}{6(n-1)}}\}.$ Hence by choosing a finite cover from $\{S_{\!\phi}(y,R)\times (-\dz^{\f{1}{n-1}},0]\}_{y\in F_\dz}$ and choosing $c_1$ such that $\f{c_1}{c_2}\leq k^{-2}$, we obtain
$$\inf_{F_\dz\times (-\dz^{\f{1}{n-1}},0]}u_1\geq c_2.$$

Fix $t_0\in(-\dz^{\f{1}{n-1}}/2,0]$. Consider the function
$$\underline{w}_{t_0}(x,t):=x_n-\phi(x)+2\l(\dz^{\f{1}{n-1}}|x'|^2+\f{\Lambda^n}{\lambda^{n-1}\dz}x_n^2\r)+(t-t_0).$$
Denote $D:=\{x_n\leq 2\dz\}\cap U$. Using similar arguments as in the proof of \cite[Lemma 5.2]{S1} and noting that on $\l(D\cap\{|x'|\leq\dz^{\f{1}{6(n-1)}}\}\r)\times \{-\dz^{\f{1}{n-1}}\}$,
$$\underline{w}_{t_0}(x,-\dz^{\f{1}{n-1}})\leq 2\dz+2\dz^{\f{1}{n-1}}\dz^{\f{1}{3(n-1)}}+2\f{\Lambda^n}{\lambda^{n-1}}4\dz-\f{1}{2}\dz^{\f{1}{n-1}}\leq 0$$
if $\dz>0$ small, we conclude that $\f{u_1}{c_2}\geq\underline{w}_{t_0}\;\mathrm{in}\;D\times [-\dz^{\f{1}{n-1}},t_0]$. In particular,
\begin{eqnarray*}
u_1(0,x_n,t_0)\geq c_2\,\underline{w}_{t_0}(0,x_n,t_0)\geq \f{c_2}{2}x_n,\quad\quad\forall x_n\in [0,c'],
\end{eqnarray*}
where $c'$ is universal.

We can apply similar arguments to any $x_0\in\partial U\cap B_c$ and $t_0\in[-\dz^{\f{1}{n-1}}/2,0]$. Hence
\begin{eqnarray*}
u_1&\geq&\f{c_2}{2}d_G,\quad\quad\mathrm{in}\;\l(\overline{U}\cap B_{c'}\r)\times [-\dz^{\f{1}{n-1}}/2,0],
\end{eqnarray*}
which implies
$$u\geq a'\,d_G\quad\mathrm{in}\;Q_\tz,\quad\quad a':=a+\f{b-a}{2}c_2.$$
\end{proof}
\end{lem}

Let $u$ be a continuous solution of \eqref{us}. Assume $t_0\in (0,T]$, $h>0$ is small and $t_0-h\geq 0$. Let $A_h$ be the sliding in the Localization Theorem and $\phi_h$ be the rescaled function of $\phi$ defined as before, we define $\mathcal{I}_{\!h} (x,t):=(h^{-\f{1}{2}}A_h x,h^{-1}(t-t_0))$ and
$$u_h(x,t):=\f{u(\mathcal{I}_{\!h}^{-1}(x,t))}{h^{1/2}}=\f{u(h^{\f{1}{2}}A_h^{-1}x,h t+t_0)}{h^{\f{1}{2}}},\quad (x,t)\in\mathcal{I}_h(Q_{\phi}((0,t_0),h)).$$
Denote $Q_{h}(\phi):=Q_{\phi}((0,0),h)$, and then
$$\mathcal{I}_{\!h}(Q_{\phi}((0,t_0),h))=Q_{1}(\phi_h)=S_{\!\phi_h}(0,1)\times (-1,0],$$
$$\mathcal{L}_{\phi_h}u_h=f_h,\quad\quad f_h(x,t):=h^{\f{1}{2}}f(h^{1/2}A_h^{-1}x,h t+t_0),\quad (x,t)\in Q_{1}(\phi_h).$$

Using this parabolic rescaling and Lemma \ref{s5:lem1}, \ref{s5:lem3}, and similar arguments as in the proof of \cite[Theorem 2.1]{S1}, we can obtain regularity of $u$ on the side of $\Omega$.

\begin{thm}\label{s5:thm2}
Under the assumptions in Lemma \ref{s5:lem1}, for any $t_0\in [2c_*,T]$ we have
\begin{eqnarray*}
|u(x,t)-D_n u(0,t_0)x_n|\leq C[|x|+|t-t_0|^{\f{1}{2}}]^{1+\alpha_1},\quad\forall (x,t)\in Q_{\phi}((0,t_0),h_1)
\end{eqnarray*}
and
\begin{eqnarray*}
|D_n u(0,t_0)|\leq C,
\end{eqnarray*}
where $\az_1\in(0,1)$ is universal and $C>0,h_1\in(0,1)$ depend only on $n,\lambda,\Lambda,\rz,c_*$.
\end{thm}

Next we assume the global information on $\Omega$ and $\phi$:

Assume 
\begin{equation}\label{5.1}
\Omega\subset B_{1/\rz}\;\mathrm{contains\;an\;interior\;ball\;of\;radius}\;\rz\;\mathrm{tangent\;to}\;\partial\Omega\;\mathrm{at\;each\;point\;on}\;\partial\Omega. 
\end{equation}
Let $\phi:\overline{\Omega}\rightarrow\mathbb{R}$, $\phi\in C^{0,1}(\overline{\Omega})\cap C^2(\Omega)$ be a convex function satisfying
\begin{equation}\label{5.2}
\mathrm{det} D^2\phi=g,\quad\quad 0<\lambda\leq g\leq\Lambda\quad\mathrm{in}\;\Omega.
\end{equation}
Assume further that $\phi$ separates quadratically from its tangent planes on $\partial\Omega$, namely,
\begin{equation}\label{5.4}
\rz|x-x_0|^2\leq\phi(x)-\phi(x_0)-\nabla\phi(x_0)\cdot (x-x_0)\leq\rz^{-1}|x-x_0|^2,\quad\quad\forall x,x_0\in\partial\Omega.
\end{equation}

Theorem \ref{s5:thm2} easily implies the estimate for general Dirichlet boundary data. 

\begin{thm}\label{s5:thm3}
Assume \eqref{5.1}-\eqref{5.4} hold. Assume $\varphi\in C^{2,1}(\overline{D})$ with $D:=\Omega\times (0,T]$. Let $u:\overline{\Omega}\times [\f{c_*}{2},T]\rightarrow\mathbb{R}$ be a continuous solution of
$$\left\{
\begin{array}{rcl}
\L u=f&&{\mathrm{in}\;\Omega\times (\f{c_*}{2},T]},\\
u=\varphi&&{\mathrm{on}\;\partial\Omega\times [\f{c_*}{2},T]}
\end{array}\right.$$
with
$$\|u\|_{\Lfz\l(\Omega\times (\f{c_*}{2},T]\r)}+\|f/\mathrm{tr}\,\Phi\|_{\Lfz\l(\Omega\times (\f{c_*}{2},T]\r)}+\|\varphi\|_{C^{2,1}(\overline{D})}\le 1.$$
Then for any $z_0=(x_0,t_0)\in\partial\Omega\times [2c_*,T]$, there exists is a linear function $l_{z_0}$ such that 
\begin{eqnarray*}
|u(x,t)-l_{z_0}(x)|\leq C[|x-x_0|+|t-t_0|^{\f{1}{2}}]^{1+\alpha_1},\quad\forall(x,t)\in Q_{\phi}(z_0,h_1)
\end{eqnarray*}
and
\begin{eqnarray*}
|D l_{z_0}|\leq C,
\end{eqnarray*}
where $\az_1\in(0,1)$ is universal and $C>0,h_1\in(0,1)$ depend only on $n,\lambda,\Lambda,\rz,c_*$.
\end{thm}

\section{Regularity near the corner of $\Omega$}\label{s6}

In this section we establish regularity for solutions $u:\l(B_\rz\cap\overline{\Omega}\r)\times [0,\rz^2]\rightarrow\mathbb{R}$ of 
\begin{equation}\label{uc}
\left\{
\begin{array}{rcl}
\L u=f&&{\mathrm{in}\;(B_\rz\cap\Omega)\times (0,\rz^2]},\\
u=\varphi&&{\mathrm{on}\;\Sigma:=\l((B_\rz\cap\partial\Omega)\times (0,\rz^2]\r)\cup\l((B_\rz\cap\overline{\Omega})\times \{0\}\r)},
\end{array}\right.
\end{equation}
where $\Omega$ and $\phi$ satisfy \eqref{s5:lem11}-\eqref{s5:lem14}. 

We first give a pointwise $C^{1+\az,\f{1+\az}{2}}$ estimate at the corner of $\Omega$. The idea is similar to that in Lemma \ref{s4:lem1}, using the good geometry of boundary sections of $\phi$ given by the Localization Theorem \ref{s5:thm1}.

\begin{lem}\label{s6:lem1}
Assume \eqref{s5:lem11}-\eqref{s5:lem14} hold. Let $u$ be a continuous solution of \eqref{uc} with
$$\|u\|_{\Lfz(\l(B_{\rz}\cap\Omega\r)\times (0, \rz^2])}+\|f/\mathrm{tr}\Phi\|_{\Lfz(\l(B_{\rz}\cap\Omega\r)\times (0,\rz^2])}\leq 1,$$
where there is a linear function $l_{(0,0)}$ such that
\begin{equation*}\label{s6:lem11}
|\varphi(x,t)-l_{(0,0)}(x)|\leq |x|^{2}+t,\quad\quad\forall (x,t)\in\Sigma,
\end{equation*}
and
\begin{equation*}\label{s6:lem11'}
|D l_{(0,0)}|\leq 1.
\end{equation*}
Then for any $0<\alpha<1$, there is a constant $c_0>0$ such that 
\begin{eqnarray*}
|u(x,t)-l_{(0,0)}(x)|&\leq& C(|x|^2+t)^{\f{1+\alpha}{2}},\quad\forall(x,t)\in \l(B_{\rz}\cap\Omega\r)\times [0,\rz^2],
\end{eqnarray*}
where $C>0,c_0>0$ depend only on $n,\lambda,\Lambda,\alpha,\rz$.

\begin{proof}
Assume $\phi(0)=0,\nabla\phi(0)=0$. For any $0<\dz<1$, the Localization Theorem \ref{s5:thm1} gives that (see \eqref{s5:lem14'}) for any $0<\dz<1$,
\begin{equation}\label{s6:lem12}
|x|^{2+\dz}\leq\phi(x)\leq |x|^{2-\dz},
\end{equation}
holds for any $x\in B_{c_0}\cap\Omega$, where $c_0$ depends only on $n,\lambda,\Lambda,\rz,\dz$.

Let $v(x,t):=\phi(x)+\mu t$ with $\mu\geq n+1/\lambda$. We have
\begin{eqnarray}\label{s6:lem13}
\L v=-\mu\mathrm{det} D^2\phi+n\mathrm{det} D^2\phi\leq -\lambda(\mu-n)\leq -1.
\end{eqnarray}
Moreover, the estimate \eqref{s6:lem12} implies that
\begin{eqnarray}\label{s6:lem14}
v(x,0)=\phi(x)\geq |x|^{2+\dz}\quad\quad x\in B_{c_0}\cap\Omega.
\end{eqnarray}

For any $\epsilon>0$, by straightforward computation we have 
$$K_1\epsilon^{-\f{\dz}{2}}\ge r^{-\dz}-\epsilon\,r^{-(2+\dz)},\quad\quad\forall r>0,$$
where $K_1:=\l[\f{\dz}{2+\dz}\r]^{1+\f{\dz}{2}}\f{2}{\dz}$. For any $x\in B_{c_0}\cap\Omega$, we have by \eqref{s6:lem14}
\begin{eqnarray}\label{s6:lem15}
|x|^2\leq\epsilon+K_1\epsilon^{-\f{\dz}{2}}|x|^{2+\dz}\leq\epsilon+K_1\epsilon^{-\f{\dz}{2}}v(x,0).
\end{eqnarray}

Define
$$\omega^{\pm}(x,t):=2[\epsilon+K_1\epsilon^{-\f{\dz}{2}}v(x,t)]+t-\f{1}{2}|x|^2\pm[u(x,t)-l(x)],\quad\quad x\in\overline{\Omega}\cap B_\rz,t\in [0,\rz^2],$$
where $l(x):=l_{(0,0)}(x)$. By \eqref{s6:lem13} we have
$$\L\omega^{\pm}=2K_1\epsilon^{-\f{\dz}{2}}\L v-\mathrm{det} D^2\phi-\mathrm{tr}\Phi[1\mp f/\mathrm{tr}\Phi]\leq
0\quad\quad\mathrm{in}\;\l(B_{c_0}\cap\Omega\r)\times (0,\rz^2].$$
On $\l(B_{c_0}\cap\Omega\r)\times\{0\}$, the estimate \eqref{s6:lem15} implies that
$$\omega^{\pm}(x,0)=2[\epsilon+K_1\epsilon^{-\f{\dz}{2}}v(x,0)]-\f{1}{2}|x|^2\pm[\varphi(x,0)-l(x)]\geq 0.$$
On $\l(\partial\Omega\cap B_{c_0}\r)\times (0,\rz^2]$, using \eqref{s6:lem15} again we have
\begin{eqnarray*}
\omega^{\pm}(x,t)&\geq&2[\epsilon+K_1\epsilon^{-\f{\dz}{2}}v(x,0)]+t-\f{1}{2}|x|^2\pm[\varphi(x,t)-l(x)]\\
&\geq&2[\epsilon+K_1\epsilon^{-\f{\dz}{2}}v(x,0)]+t-\f{1}{2}|x|^2-[|x|^{2}+t]\geq 0.
\end{eqnarray*}
On $\l(\partial B_{c_0}\cap\Omega\r)\times (0,\rz^2]$, we use  \eqref{s6:lem14} and find
\begin{eqnarray*}
\omega^{\pm}(x,t)&\geq&2K_1\epsilon^{-\f{\dz}{2}}|x|^{2+\dz}-\f{1}{2}|x|^2-[|\varphi(x,t)|+|l(0)|+|Dl||x|]\\
&\geq&2K_1\epsilon^{-\f{\dz}{2}}c_0^{2+\dz}-\f{1}{2}c_0^2-3\geq 0,
\end{eqnarray*}
where we choose $\ez>0$ such that $K_1\epsilon^{-\f{\dz}{2}}\geq 2c_0^{-(2+\dz)}$. Then the maximum principle implies that
$$\omega^{\pm}\geq 0\quad\quad\mathrm{in}\;\l(B_{c_0}\cap\Omega\r)\times [0,\rz^2].$$
Hence for any $\ez$ such that $K_1\epsilon^{-\f{\dz}{2}}\geq 2c_0^{-(2+\dz)}$ we have
\begin{eqnarray}\label{s6:lem16}
|u(x,t)-l(x)|\leq2[\epsilon+K_1\epsilon^{-\f{\dz}{2}}v(x,t)]+t,\quad\quad\forall (x,t)\in \l(B_{c_0}\cap\Omega\r)\times [0,\rz^2],
\end{eqnarray}
which gives
\begin{eqnarray*}
|u(x,t)-l(x)|\leq C[v(x,t)]^{\f{2}{2+\dz}}+t\leq C[|x|^{\f{2(2-\dz)}{2+\dz}}+t^{\f{2}{2+\dz}}],
\end{eqnarray*}
where $C$ depends only on $n,\lz,\Lz,\rz,\dz$. For any $\alpha<1$, we choose $\dz>0$ small such that $\f{2(2-\dz)}{2+\dz}>1+\alpha$ and $\f{2}{2+\dz}>\f{1+\alpha}{2}$. The result follows.
\end{proof}
\end{lem}

Lemma \ref{s6:lem1}, Theorem \ref{s5:thm2} and a localization process give estimates of $u$ near the corner of $\Omega$.

\begin{lem}\label{s6:lem2}
Assume \eqref{s5:lem11}-\eqref{s5:lem14} hold. Let $u$ be a continuous solution of \eqref{uc} with $\varphi=0$ and
$$\|u\|_{L^\infty((B_\rz\cap\Omega)\times (0,\rz^2])}+\|f/\mathrm{tr}\Phi\|_{L^\infty((B_\rz\cap\Omega)\times (0,\rz^2])}\leq 1.$$
Then for any $t_0\leq h_0$ we have
\begin{eqnarray*}
|u(x,t)-D_n u(0,t_0)x_n|\leq C[|x|+|t-t_0|^{\f{1}{2}}]^{1+\alpha_1},\quad\forall(x,t)\in Q_{\phi}((0,t_0),h_1t_0), 
\end{eqnarray*}
and
\begin{eqnarray*}
|D_n u(0,t_0)|\leq Ct_0^{\f{\alpha_1}{2}},
\end{eqnarray*}
where $C>0,h_0,h_1,\alpha_1\in(0,1)$ are universal.

\begin{proof}
Recall the rescaled function $\phi_h$ and section $S_{\!\phi_h}(0,1)$ in Section \ref{s5}. For any $h\leq h_0=h_0(n,\lambda,\Lambda,\rz)$, $\phi_h$ satisfies in $S_{\!\phi_h}(0,1)$ the hypotheses of the Localization Theorem \ref{s5:thm1} at all points on $\partial S_{\!\phi_h}(0,1)\cap B_{c_0}$ for some $c_0$ universal. Define $I_h(x,t):=(h^{-\f{1}{2}}A_h x,h^{-1}t)$ and
$$Q'_h(\phi):=S_{\!\phi}(0,h)\times (0,h],$$
then
$$\mathcal{I}_hQ'_h(\phi)=Q'_1(\phi_h)=S_{\!\phi_h}(0,1)\times (0,1].$$
Define
$$u_h(y,s):=\f{u(h^{\f{1}{2}}A_h^{-1}y,hs)}{h^{\f{1}{2}}},$$
then
$$\mathcal{L}_{\phi_h}u_h(y,s)=f_h(y,s):=h^{\f{1}{2}}f(h^{\f{1}{2}}A_h^{-1}y,hs)\quad\quad\mathrm{in}\;Q'_1(\phi_h).$$
Note
$$S_{\!\phi_h}(0,1)\supset B_{k}\cap\overline{\Omega_h}$$
with $k$ universal and
$$u_h=0\quad\quad\mathrm{on}\;\l((\partial S_{\!\phi_h}(0,1)\cap B_{k})\times (0,1]\r)\cup\l((S_{\!\phi_h}(0,1)\cap B_{k})\times\{0\}\r).$$
Applying Theorem \ref{s5:thm2}, for any $s_0\in[\f{1}{2},1]$ and any $(y,s)\in Q_{\phi_h}((0,s_0),h_1)$, we have
\begin{eqnarray*}\label{s6:rem21}
|u_h(y,s)-D_n u_h(0,s_0)y_n|&\leq&C[\|u_h\|_{\Lfz((S_{\!\phi_h}(0,1)\cap B_k)\times (0,1])}+\|f_h/\mathrm{tr}\Phi_h\|_{\Lfz((S_{\!\phi_h}(0,1)\cap B_k)\times (0,1])}]\nonumber\\
&\cdot&[|y|+|s-s_0|^{\f{1}{2}}]^{1+\alpha_1},
\end{eqnarray*}
and
\begin{eqnarray*}
|D_n u_h(0,s_0)|&\leq&C[\|u_h\|_{\Lfz((S_{\!\phi_h}(0,1)\cap B_k)\times (0,1])}+\|f_h/\mathrm{tr}\Phi_h\|_{\Lfz((S_{\!\phi_h}(0,1)\cap B_k)\times (0,1])}],
\end{eqnarray*}
where $C>0,0<h_1,\alpha_1<1$ are universal.

By Lemma \ref{s6:lem1}, for any $\alpha\in(0,1)$ and $(x,t)\in Q'_h(\phi)$, we obtain
\begin{eqnarray}\label{s6:lem21}
|u(x,t)|\leq C\l(|x|^2+t\r)^{\f{1+\alpha}{2}}\leq Ch^{\f{1+\alpha}{2}}|\log h|^{1+\alpha},
\end{eqnarray}
where $C=C(n,\lz,\Lz,\rz,\az)>0$. Moreover,
\begin{eqnarray}\label{s6:lem22}
\|f_h/\mathrm{tr}\Phi_h\|_{\Lfz\l(\l(B_{k}\cap S_{\!\phi_h}(0,1)\r)\times (0,1]\r)}
\leq Ch^{\f{1}{2}}|\log h|^{2}\|f/\mathrm{tr}\Phi\|_{L^\infty(Q'_h(\phi))}\leq Ch^{\f{1}{2}}|\log h|^{2}
\end{eqnarray}
for some constant $C=C(n,\lz,\Lz,\rz)>0$.

For any $t_0\in [\f{h}{2},h]$, set
$$l_{(0,t_0)}(x):=h^{\f{1}{2}}D_n u_h(0,h^{-1}t_0)(h^{-\f{1}{2}}A_h x)\cdot e_n=D_n u_h(0,h^{-1}t_0)x_n.$$
For any $(x,t)\in Q_{\phi}((0,t_0),h_1h)$, we have by \eqref{s6:lem21} and \eqref{s6:lem22} that
\begin{eqnarray*}
|u(x,t)-l_{(0,t_0)}(x)|&\leq&Ch^{\f{1}{2}}\l[h^{\f{\alpha}{2}}|\log h|^{1+\alpha}
+h^{\f{1}{2}}|\log h|^{2}\r]\cdot[|h^{-\f{1}{2}}A_h x|+|h^{-1}(t-t_0)|^{\f{1}{2}}]^{1+\alpha_1}\\
&\leq&Ch^{\f{1}{2}}h^{\f{\alpha}{2}}|\log h|^{2}h^{-\f{1+\alpha_1}{2}}|\log h|^{1+\alpha_1}[|x|+|t-t_0|^{\f{1}{2}}]^{1+\alpha_1}\\
&\leq&C[|x|+|t-t_0|^{\f{1}{2}}]^{1+\alpha_1},
\end{eqnarray*}
and
\begin{eqnarray*}
|D_n u_h(0,s_0)|&\leq&C\l[h^{\f{\alpha}{2}}|\log h|^{1+\alpha}
+h^{\f{1}{2}}|\log h|^{2}\r]\le Ch^{\f{\az_1}{2}},
\end{eqnarray*}
where $C=C(n,\lz,\Lz,\rz)$ and we take $\az>\alpha_1$. Moreover, from the first estimate we conclude that $D_n u(0,t_0)=D_n u_h(0,h^{-1}t_0)$.

For any $t_0\leq h_0$, take $h:=t_0\leq h_0$, then the last two estimates imply the desired conclusion.
\end{proof}
\end{lem}

Lemma \ref{s6:lem2} implies the estimate for general Dirichlet boundary data.

\begin{thm}\label{s6:thm1}
Assume \eqref{5.1}-\eqref{5.4} hold. Denote $D:=\Omega\times(0,T]$. Let $\varphi\in C^{2,1}(\overline{D})$ and $u$ be a continuous solution of
$$\left\{
\begin{array}{rcl}
\L u=f&&\mathrm{in}\;D,\\
u=\varphi&&\mathrm{on}\;\partial_pD
\end{array}\right.$$
with
$$\|u\|_{\Lfz(D)}+\|f/\mathrm{tr}\,\Phi\|_{\Lfz(D)}+\|\varphi\|_{C^{2,1}(\overline{D})}\le 1.$$
Then, for any $z_0=(x_0,t_0)\in\partial\Omega\times [0,h_0]$, there exists is a linear function $l_{z_0}$ such that 
\begin{eqnarray*}
|u(x,t)-l_{z_0}(x)|\leq C[|x-x_0|+|t-t_0|^{\f{1}{2}}]^{1+\alpha_1},\quad\forall(x,t)\in Q_{\phi}(z_0,h_1t_0)
\end{eqnarray*}
and
\begin{eqnarray*}
\sup_{(x,t)\in\partial\Omega\times [0,h_0]}|D l_{(x,t)}|+\sup_{\substack{(x,t),(y,s)\in\partial\Omega\times[0,h_0]\\ (x,t)\neq(y,s)}}\f{|D l_{(x,t)}-D l_{(y,s)}|}{(|x-y|+|t-s|^{\f{1}{2}})^{\alpha_1}}\le C,
\end{eqnarray*}
where $C>0,h_0,h_1,\alpha_1\in(0,1)$ are universal.
\end{thm}

\section{Proof of Theorem 2}\label{s7}

In this section we give the complete proof of Theorem $2$. We denote by $c,C$ constants depending only on $n,\lz,\Lz,\rz,T$ and the modulus of continuity of $g$. Their values may change from line to line whenever there is no confusion.

We may assume
$$\|u\|_{\Lfz(D)}+\|f\|_{\Lfz(D)}+\|\varphi\|_{C^{2,1}(\overline{D})}\leq 1.$$

Let $y\in\Omega$. Assume $\Omega\subset\mathbb{R}^n_+, 0\in\partial\Omega$ and the maximal interior section $S_{\!\phi}(y,\bar{h})$ is tangent to $\partial\Omega$ at $0$. By Proposition \ref{s5:prop1}, we have
$$\mathrm{dist}(y,\partial\Omega)\sim\bar{h}^{\f{1}{2}},\quad\quad |\nabla\phi(y)|\sim\bar{h}^{\f{1}{2}},$$
$$c E\subset S_{\!\phi}(y,\bar{h})-y\subset C E,$$
where
$$E:=\bar{h}^{\f{1}{2}}A_{\bar{h}}^{-1}B_1,\quad\quad \|A_{\bar{h}}\|,\|A_{\bar{h}}^{-1}\|\leq C\log|\bar{h}|.$$
Define $Tx:=\bar{h}^{-\f{1}{2}}A_{\bar{h}}(x-y)$ and $T_{\!p}(x,t):=(Tx,\bar{h}^{-1}t)$,
\begin{eqnarray*}
\tilde{\phi}(\tilde{x}):&=&\f{1}{\bar{h}}[\phi(T^{-1}\tilde{x})-l_{\phi,y} (T^{-1}\tilde{x})-\bar{h}],\\
\tilde{u}(\tilde{x},\tilde{t}):&=&\f{1}{\bar{h}^{\f{1}{2}}}[u(T^{-1}\tilde{x},\bar{h}\tilde{t})-l_{(0,0)}(T^{-1}\tilde{x})],
\end{eqnarray*}
where $l_{(0,0)}$ is from Lemma \ref{s6:lem1}. Let $K_0>1$ be a constant to be chosen later. Then
\begin{eqnarray*}
&&B_c\times (0,K_0]\subset T_{\!p}(S_{\!\phi}(y,\bar{h})\times (0,K_0\bar{h}])=S_{\!\tilde{\phi}}(0,1)\times (0,K_0]\subset B_C\times (0,K_0],\\
&&\mathcal{L}_{\tilde{\phi}}\tilde{u}(\tilde{x},\tilde{t})=\tilde{f}(\tilde{x},\tilde{t})
:=\bar{h}^{\f{1}{2}}f(T^{-1}\tilde{x},\bar{h}\tilde{t})\quad\quad\mathrm{in}\;S_{\!\tilde{\phi}}(0,1)\times (0,K_0].
\end{eqnarray*}
From Lemma \ref{s6:lem1}, for any $\alpha\in(0,1)$ we have
\begin{eqnarray}\label{s5:thm414'}
|u(x,t)-l_{(0,0)}(x)|\leq C[|x|+|t|^{\f{1}{2}}]^{1+\alpha},\quad\forall (x,t)\in D,
\end{eqnarray}
Define $\tilde{\varphi}:=\tilde{u}(\cdot,0)$, then by straightforward computation we obtain
\begin{eqnarray}\label{s5:thm414''}
\|\tilde{\varphi}\|_{C^{1,1}(\overline{S_{\!\tilde{\phi}}(0,1)})}\leq C\bar{h}^{\f{1}{2}}|\log\bar{h}|^{2}.
\end{eqnarray}
Theorem \ref{s4:thm3} gives that $\tilde{u}\in C^{1+\alpha',\f{1+\alpha'}{2}}\l(\overline{S_{\!\tilde{\phi}}(0,\f{1}{2})}\times [0,K_0]\r)$ for any $\alpha'<1$ and
\begin{eqnarray*}
\|\tilde{u}\|_{C^{1+\alpha',\f{1+\alpha'}{2}}\l(\overline{S_{\!\tilde{\phi}}(0,\f{1}{2})}\times \l[0,K_0\r]\r)}\leq C[\|\tilde{u}\|_{\Lfz(S_{\!\tilde{\phi}}(0,1)\times (0,K_0])}+\|\tilde{f}\|_{\Lfz(S_{\!\tilde{\phi}}(0,1)\times (0,K_0])}+\|\tilde{\varphi}\|_{C^{1,1}(\overline{S_{\!\tilde{\phi}}(0,1)})}]\leq\bar{h}^{\f{\alpha}{4}},
\end{eqnarray*}
where we use \eqref{s5:thm414'} and \eqref{s5:thm414''}. Hence,
\begin{eqnarray*}
\|D u\|_{L^\infty\l(\overline{S_{\!\phi}\l(y,\f{\bar{h}}{2}\r)}\times(0,K_0\bar{h}]\r)}&\leq&|D l_{(0,0)}|
+\|A_{\bar{h}}\|\,\|D\tilde{u}\|_{L^\infty\l(\overline{S_{\!\tilde{\phi}}(0,\f{1}{2})}\times \l[0,K_0\r]\r)}\leq C,
\end{eqnarray*}
and for any $(x,t_1),(x,t_2)\in S_{\!\phi}\l(y,\f{\bar{h}}{2}\r)\times(0,K_0\bar{h}]$,
\begin{eqnarray*}
|u(x,t_1)-u(x,t_2)|&=&\bar{h}^{\f{1}{2}}|\tilde{u}(Tx,\bar{h}^{-1}t_1)-\tilde{u}(Tx,\bar{h}^{-1}t_2)|\\
&\leq&\bar{h}^{\f{1}{2}}\bar{h}^{\f{\alpha}{4}}\bar{h}^{-\f{1+\alpha'}{2}}|t_1-t_2|^{\f{1+\alpha'}{2}}\leq |t_1-t_2|^{\f{1+\alpha'}{2}},
\end{eqnarray*}
by taking $\alpha'<\f{\alpha}{2}$. For any $(x_i,t_i)\in S_{\!\phi}\l(y,\f{\bar{h}}{2}\r)\times(0,K_0\bar{h}],\;i=1,2$,
\begin{eqnarray*}
|D u(x_1,t_1)-D u(x_2,t_2)|
&=&|(A_{\bar{h}})^t[D\tilde{u}(Tx_1,\bar{h}^{-1}t_1)-D\tilde{u}(Tx_2,\bar{h}^{-1}t_2)]|\\
&\leq&C|\log\bar{h}|\,\bar{h}^{\f{\alpha}{4}}\,[|\bar{h}^{-\f{1}{2}}A_{\bar{h}}(x_1-x_2)|+\bar{h}^{-\f{1}{2}}|t_1-t_2|^{\f{1}{2}}]^{\alpha'}\\
&\leq&C\bar{h}^{\f{\alpha}{4}-\f{\alpha'}{2}}|\log\bar{h}|^{1+\alpha'}[|x_1-x_2|+|t_1-t_2|^{\f{1}{2}}]^{\alpha'}\\
&\leq&[|x_1-x_2|+|t_1-t_2|^{\f{1}{2}}]^{\alpha'}.
\end{eqnarray*}
Hence we obtain
\begin{eqnarray}\label{s5:thm415'}
\|u\|_{C^{1+\alpha',\f{1+\alpha'}{2}}\l(\overline{S_{\!\phi}(y,\f{\bar{h}}{2})}\times \l[0,K_0\bar{h}\r]\r)}\leq C.
\end{eqnarray}
Denote $l_{(y,t)}:=l_{u,(y,t)}$ for $t\in(0,T]$. For any $t_0\in (0,K_0\bar{h}]$ and $(x,t)\in\overline{S_{\!\phi}(y,\f{\bar{h}}{2})}\times[0,K_0\bar{h}]$,
\begin{eqnarray}\label{s5:thm416}
|u(x,t)-l_{(y,t_0)}(x)|&=&|u(x,t)-u(y,t_0)-D u(y,t_0)\cdot (x-y)|\nonumber\\
&\leq&C[|x-y|+|t-t_0|^{\f{1}{2}}]^{1+\alpha'}.
\end{eqnarray}
Then by \eqref{s5:thm414'} and \eqref{s5:thm416},
\begin{eqnarray*}
|D l_{(0,0)}-D u(y,t_0)|\leq C |y|^{\f{\alpha'}{2}}.
\end{eqnarray*}
Let $\az_1\in(0,1)$ be the small constant given by Theorem \ref{s6:thm1}, then 
\begin{eqnarray*}
|D u(y,t_0)-D l_{(0,t_0)}|&\leq&|D l_{(0,0)}-D u(y,t_0)|+|D l_{(0,0)}-D l_{(0,t_0)}|\leq C|y|^{\f{\az'}{2}},
\end{eqnarray*}
where we use Theorem \ref{s6:thm1} and take $\f{\az'}{2}<\az_1$. Combining these estimates we find
\begin{eqnarray}\label{s5:thm417}
\sup_{t\in[0,K_0\bar{h}]}\f{|D u(y,t)-D l_{(0,t)}|}{|y|^{\f{\az'}{2}}}\leq C.
\end{eqnarray}

Let $h_0,h_1$ be the small constants given by Theorem \ref{s6:thm1}. Fix $t_0\in (K_0\bar{h},h_0]$. Define
$$v(\tilde{x},\tilde{t}):=\f{1}{\bar{h}^{\f{1}{2}}}\l[u-l_{(0,t_0)}\r](\bar{h}^{\f{1}{2}}A_{\bar{h}}^{-1}\tilde{x}+y,\bar{h}\tilde{t}+t_0),$$
then
$$\mathcal{L}_{\tilde{\phi}}v(\tilde{x},\tilde{t})=F(\tilde{x},\tilde{t}):=\bar{h}^{\f{1}{2}}f(\bar{h}^{\f{1}{2}}A_{\bar{h}}^{-1}\tilde{x},\bar{h}\tilde{t}+t_0)\quad\quad\mathrm{in}\;Q_{\tilde{\phi}}((0,0),1).$$
The interior estimate Theorem $1$ shows for $\alpha'<1$ as above,
\begin{eqnarray*}
\|v\|_{C^{1+\alpha',\f{1+\alpha'}{2}}\l(\overline{Q_{\tilde{\phi}}((0,0),\f{1}{2})}\r)}\leq C[\|v\|_{\Lfz(Q_{\tilde{\phi}}((0,0),1))}+\|F\|_{\Lfz(Q_{\tilde{\phi}}((0,0),1))}].
\end{eqnarray*}
Applying Theorem \ref{s6:thm1}, there exists is a linear function $l_{(0,t_0)}$ such that 
\begin{eqnarray}\label{s5:thm418}
|u(x,t)-l_{(0,t_0)}(x)|\leq C[|x|+|t-t_0|^{\f{1}{2}}]^{1+\alpha_1},\quad\forall (x,t)\in Q_{\phi}((0,t_0),h_1t_0).
\end{eqnarray}
By \cite[Proposition 4.1]{N},
$$S_{\!\phi}(y,\bar{h})\subset S_{\!\phi}(0,\tz_*\bar{h}),$$
where $\tz_*$ depends only on $n,\lz,\Lz,\rz$.

Now we take $K_0:=\tz_*h_1^{-1}$, then $h_1t_0>h_1K_0\bar{h}=\tz_*\bar{h}$, we have by \eqref{s5:thm418}
\begin{eqnarray}\label{s5:thm423}
|u(x,t)-l_{(0,t_0)}(x)|\leq C[|x|+|t-t_0|^{\f{1}{2}}]^{1+\alpha_1},\quad\forall(x,t)\in Q_{\phi}((y,t_0),\bar{h}).
\end{eqnarray}
Hence, for any $(\tilde{x},\tilde{t})\in Q_{\tilde{\phi}}((0,0),1)$,
\begin{eqnarray*}
|v(\tilde{x},\tilde{t})|
&\leq&C\bar{h}^{-\f{1}{2}}[|\bar{h}^{\f{1}{2}}A_{\bar{h}}^{-1}\tilde{x}+y|+\bar{h}^{\f{1}{2}}|\tilde{t}|^{\f{1}{2}}]^{1+\alpha_1}
\leq\bar{h}^{\f{\alpha_1}{4}}.
\end{eqnarray*}
Similar to the proof of \eqref{s5:thm415'} we have
\begin{eqnarray}\label{s5:thm420}
\|u\|_{C^{1+\alpha',\f{1+\alpha'}{2}}\l(\overline{Q_{\phi}((y,t_0),\f{\bar{h}}{2})}\r)}\leq C
\end{eqnarray}
by choosing $\alpha'<\f{\alpha_1}{2}$. For any $(x,t)\in Q_{\phi}((y,t_0),\f{\bar{h}}{2})$, we have from \eqref{s5:thm420} that
\begin{eqnarray}\label{s5:thm422}
|u(x,t)-l_{(y,t_0)}(x)|\leq C[|x-y|+|t-t_0|^{\f{1}{2}}]^{1+\alpha'}.
\end{eqnarray}
This together with \eqref{s5:thm423} implies
\begin{eqnarray*}
|D l_{(x_0,t_0)}-D u(y,t_0)|\leq C |y|^{\f{\alpha'}{2}}.
\end{eqnarray*}
Hence,
\begin{eqnarray}\label{s5:thm424}
\sup_{t\in[K_0\bar{h},h_0]}\f{|D u(y,t)-D l_{(0,t)}|}{|y|^{\f{\alpha'}{2}}}\leq C.
\end{eqnarray}

Note that \eqref{s5:thm418} also holds if $t_0\ge h_0$ (by applying Theorem \ref{s5:thm3}), hence \eqref{s5:thm420} and \eqref{s5:thm424} also hold for $t_0\ge h_0$.

Let $(x,t),(y,s)\in\Omega\times [0,T]$. Denote $\bar{h}_x:=\bar{h}(x),\bar{h}_y:=\bar{h}(y)$ and assume $\bar{h}_y\leq\bar{h}_x$. Let $x^*, y^*\in\partial\Omega$ such that $x^*\in\partial S_{\!\phi}(x,\bar{h}_x)\cap\partial\Omega$ and $y^*\in\partial S_{\!\phi}(y,\bar{h}_y)\cap\partial\Omega,$. We consider these cases:\\

$\mathbf{Case\;1.}$ $y\in S_{\!\phi}\l(x,\f{\bar{h}_x}{2}\r),\;|t-s|<\f{\bar{h}_x}{2},\;t,s\leq K_0\bar{h}_x$. Then
$$(x,t),(y,s),(x,s)\in S_{\!\phi}\l(x,\f{\bar{h}_x}{2}\r)\times (0,K_0\bar{h}_x].$$
Hence by \eqref{s5:thm415'},
\begin{eqnarray*}
|D u(x,t)-D u(y,s)|&\leq&C[|x-y|+|t-s|^{\f{1}{2}}]^{\alpha'},\\
|u(x,t)-u(x,s)|&\leq&C|t-s|^{\f{1+\alpha'}{2}}.
\end{eqnarray*}

$\mathbf{Case\;2.}$ $y\in S_{\!\phi}\l(x,\f{\bar{h}_x}{2}\r),\;|t-s|<\f{\bar{h}_x}{2},\;\max\{t,s\}>K_0\bar{h}_x$. For example, if $s>K_0\bar{h}_x,\;t\leq s$, then
$$(y,s),(x,t)\in Q_{\phi}\l((x,s),\f{\bar{h}_x}{2}\r).$$
Hence by \eqref{s5:thm420},
\begin{eqnarray*}
|D u(x,t)-D u(y,s)|&\leq&|D u(y,s)-D u(x,s)|+|D u(x,s)-D u(x,t)|\\
&\leq&C[|x-y|+|t-s|^{\f{1}{2}}]^{\alpha'},\\
|u(x,t)-u(x,s)|&\leq&C|t-s|^{\f{1+\alpha'}{2}}.
\end{eqnarray*}

$\mathbf{Case\;3.}$ $y\in S_{\!\phi}\l(x,\f{\bar{h}_x}{2}\r),\;|t-s|\geq\f{\bar{h}_x}{2}$.\\

For any $\dz\in(0,\f{1}{2})$ small we have
$$|y^*-y|\leq C\bar{h}_y^{\f{1}{2}-\dz}\leq C\bar{h}_x^{\f{1}{2}-\dz}\leq C|s-t|^{\f{1}{2}-\dz},$$
$$|x^*-x|\leq C\bar{h}_x^{\f{1}{2}-\dz}\leq C|s-t|^{\f{1}{2}-\dz},$$
then by \eqref{s5:thm417}, \eqref{s5:thm424} and Theorem \ref{s6:thm1},
\begin{eqnarray*}
|D u(x,t)-D u(y,s)|&\leq&|D u(x,t)-D l_{(x^*,t)}|+|D l_{(x^*,t)}-D l_{(y^*,s)}|+|D l_{(y^*,s)}-D u(y,s)|\\
&\leq&C[|x-x^*|^{\f{\alpha'}{2}}+(|x^*-y^*|+|t-s|^{\f{1}{2}})^{\alpha_1}+|y^*-y|^{\f{\alpha'}{2}}],\\
&\leq&C[|x-y|+|t-s|^{\f{1}{4}}]^{\f{\alpha'}{2}},\\
|u(x,t)-u(x,s)|&\leq&|u(x,t)-l_{(x^*,t)}(x)|+|l_{(x^*,t)}(x)-l_{(x^*,s)}(x)|+|u(x,s)-l_{(x^*,s)}(x)|\\
&=&|(D u(\xi,t)-D l_{(x^*,t)})\cdot(x-x^*)|+|(D u(\eta,s)-D l_{(x^*,s)})\cdot(x-x^*)|\\
&+&|u(x^*,t)-u(x^*,s)+(D l_{(x^*,t)}-D l_{(x^*,s)})\cdot(x-x^*)|\\
&\leq&C[|x-x^*|^{1+\f{\alpha'}{2}}+|t-s|+|t-s|^{\f{\alpha_1}{2}}|x-x^*|]\\
&\leq&C|t-s|^{(\f{1}{2}-\dz)(1+\f{\az'}{2})}\le C|t-s|^{\f{1}{2}(1+\f{\az'}{4})},
\end{eqnarray*}
where $\xi,\eta$ are some points in the segment joining $x$ and $x^*$, and we choose $\dz$ sufficiently small.\\

$\mathbf{Case\;4.}$ $y\notin S_{\!\phi}(x,\bar{h}_x/C^*)$.\\

Since $\phi\in C^{1,\bz}(\overline{\Omega})$ for some $\bz=\bz(n,\lz,\Lz,\rz)\in(0,1)$ (\cite[Proposition 2.6]{S1}), we have
$$|x-y|\geq c\bar{h}_x^{\f{1}{1+\bz}}\geq\bar{h}_y^{\f{1}{1+\bz}},$$
it follows that
$$|x^*-x|\leq C\bar{h}_x^{\f{1}{4}}\leq C|x-y|^{\f{1+\bz}{4}},$$
and
$$|y^*-y|\leq C\bar{h}_y^{\f{1}{4}}\leq C|x-y|^{\f{1+\bz}{4}},$$
then
\begin{eqnarray*}
|D u(x,t)-D u(y,s)|&\leq&|D u(x,t)-D l_{(x^*,t)}|+|D l_{(x^*,t)}-D l_{(y^*,s)}|+|D l_{(y^*,s)}-D u(y,s)|\\
&\leq&C[|x-x^*|^{\f{\alpha'}{2}}+(|x^*-y^*|+|t-s|^{\f{1}{2}})^{\alpha_1}+|y^*-y|^{\f{\alpha'}{2}}],\\
&\leq&C[|x-y|^{\f{1+\bz}{4}}+|t-s|^{\f{1}{2}}]^{\f{\alpha'}{2}}.
\end{eqnarray*}

Combining all these cases, the proof of Theorem $2$ is complete.
\qed

LMAM, School of Mathematical Sciences,
 Peking University, Beijing, 100871,
 P. R. China

 Lin Tang,\quad
 E-mail address:  tanglin@math.pku.edu.cn

Qian Zhang,\quad
E-mail address: 1401110018@pku.edu.cn
\end{document}